\documentclass[preprint,12pt]{elsarticle}

\usepackage{amssymb}
\usepackage{color}
\usepackage{tabularx}

\newcommand{\cC}{{\mathcal C}}
\newcommand{\cE}{{\mathcal E}}

\newcommand{\cN}{{\mathcal N}}

\newcommand{\cT}{{\mathcal T}}

\newtheorem{theorem}{Theorem}[section]
\newtheorem{lemma}[theorem]{Lemma}

\newtheorem{corollary}[theorem]{Corollary}
\newtheorem{observation}[theorem]{Observation}
\newtheorem{sublemma}{}[theorem]
\newproof{proof}{Proof}   

\journal{}

\begin{document}

\begin{frontmatter}

\title{Non-essential arcs in phylogenetic networks}

\author[label1]{Simone Linz}
\ead{s.linz@auckland.ac.nz}
\author[label2]{Charles Semple}
\ead{charles.semple@canterbury.ac.nz}

\address[label1]{School of Computer Science, University of Auckland, New Zealand}
\address[label2]{School of Mathematics and Statistics, University of Canterbury, New Zealand}

\begin{abstract}
In the study of rooted phylogenetic networks, analyzing the set of rooted phylogenetic trees that are embedded in such a network is a recurring task. From an algorithmic viewpoint, this analysis almost always requires an exhaustive search of a particular multiset $S$ of rooted phylogenetic trees that are embedded in a rooted phylogenetic network $\cN$. Since the size of $S$ is exponential in the number of reticulations of $\cN$, it is consequently of interest to keep this number as small as possible but without loosing any element of $S$. In this paper, we take a first step towards this goal by introducing the notion of a non-essential arc of $\cN$, which is an arc whose deletion from $\cN$ results in a rooted phylogenetic network $\cN'$ such that the sets of rooted phylogenetic trees that are embedded in $\cN$ and $\cN'$ are the same. We investigate the popular class of tree-child networks and characterize which arcs are non-essential. This characterization is based on a  family of directed graphs. Using this novel characterization, we show that identifying and deleting all non-essential arcs in a tree-child network takes time that is cubic in the number of leaves of the network. Moreover, we show that deciding if a given arc of an arbitrary phylogenetic network is non-essential is $\Pi_2^P$-complete.
\end{abstract}


\end{frontmatter}

\section{Introduction}
As a generalization of rooted phylogenetic trees to rooted directed acyclic graphs with labeled leaves, rooted phylogenetic networks provide advanced opportunities to accurately represent ancestral relationships between entities such as species, viruses, and cancer cells~\cite{gusfield14,huson11}. This is particularly true if the evolutionary history includes non-treelike processes such as hybridization, lateral gene transfer, or recombination that cannot be represented by a single phylogenetic tree. However, since the evolutionary history of a single gene or short DNA fragment is, in most cases, correctly described by a tree, the set of rooted phylogenetic trees that are embedded in a network are of recurring interest in the reconstruction and analysis of rooted phylogenetic networks. 

For example, one way of generalizing the popular tree reconstruction methods parsimony and likelihood~\cite{felsenstein04} to rooted phylogenetic networks is by scoring each phylogenetic tree that is embedded in a given network instead of scoring the network directly~\cite{hein90,jin06,nakhleh05}. Another question that considers tree embeddings is the fundamental question of whether or not a given phylogenetic tree is embedded in a given phylogenetic network. This problem, known as {\sc Tree Containment}, was first introduced by Kanj et al.~\cite{kanj08} and continues to be actively studied (e.g.~\cite{bordewich16,bentert18,gambette18,iersel18}). Computing the parsimony or likelihood score of a rooted phylogenetic network as well as solving {\sc Tree Containment} are, without imposing additional structural properties on the input, three NP-hard problems~\cite{fischer15,jin06,jin09,kanj08,nguyen07}. Other computationally hard problems that consider embeddings of rooted phylogenetic trees in a network include counting the number of phylogenetic trees that are embedded in a rooted phylogenetic network and deciding if a rooted phylogenetic tree is a so-called base tree of a rooted phylogenetic network~ \cite{anaya16,linz13}.

\begin{figure}
\center
\scalebox{1}{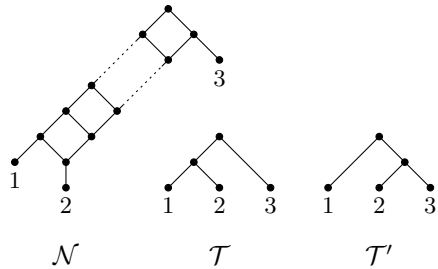}
\caption{A rooted phylogenetic network $\cN$ that, regardless of the number of reticulations, only embeds the two rooted phylogenetic trees $\cT$ and $\cT'$ shown on the right-hand side. If $\cN$ has $k$ reticulations, then there are exactly $2^k-1$ distinct combinations of arcs to delete from $\cN$ that each yield $\cT$.}
\label{fig:stack}
\end{figure}

Due to the computational complexity of the aforementioned problems, it is not surprising that algorithms to solve these problems often exhaustively search through a multiset of all rooted phylogenetic trees that are embedded in a rooted phylogenetic network $\cN$. We next describe a standard procedure to obtain this multiset of trees for when each non-root vertex in $\cN$ has in-degree one or two.  A straightforward generalization of the process can be applied to networks with vertices of higher in-degree. By deleting exactly one incoming arc for each vertex in $\cN$ that has in-degree two, and repeatedly suppressing vertices of in-degree one and out-degree one, deleting unlabeled vertices with in-degree one and out-degree zero, and deleting vertices with in-degree zero and out-degree one, one obtains a rooted phylogenetic tree that is embedded in $\cN$. Indeed, each  rooted phylogenetic tree that is embedded in $\cN$ can be obtained in this way. Hence, if $\cN$ has $k$ vertices of in-degree two, the systematic process of deleting  subsets of arcs of size $k$ in $\cN$ generates a multiset $S$ of rooted phylogenetic trees that are embedded in $\cN$. If a rooted phylogenetic tree $\cT$ of $S$ has multiplicity $c$, then this implies that there are $c$ combinations of arcs to delete in $\cN$ that each yield $\cT$.  Although the cardinality of $S$ is $2^k$, it is known that the number of elements in $S$ is often smaller than $2^k$~\cite{cordue14}.  For instance, it is easy to construct a rooted phylogenetic network that embeds exactly two (distinct) rooted phylogenetic trees for any given $k$ (see Figure~\ref{fig:stack}). In Figure~\ref{fig:stack}, as well as in all other figures of this paper, arcs are directed down the page.

Given the potential difference between the cardinality of $S$ and the number of elements in $S$, we ask the following question. Does there exist an arc $e$ in $\cN$ such that, by deleting $e$ from $\cN$, we obtain a rooted phylogenetic network $\cN'$ whose set of embedded trees is the same as the set of embedded trees of $\cN$? We make the process of obtaining $\cN'$ from $\cN$ by deleting $e$ more precise in Section~\ref{sec:preliminaries}. If there exists such an arc $e$, we call $e$ a {\it non-essential} arc of $\cN$. Furthermore, if $e$ is directed into a vertex of in-degree two in $\cN$, then, from an algorithmic viewpoint, the existence of $e$ implies that $|S|=2|S'|$, where $S'$ is the multiset of rooted phylogenetic trees that are embedded in $\cN'$. Thus, the deletion of $e$ in $\cN$ reduces the running time of an algorithm that searches $S$ by up to one half. It is consequently of interest to explore ways that simplify $\cN$ by means of arc deletions without losing any element of $S$.

In this paper, we take a first step in this exploration and characterize non-essential arcs in tree-child networks, which are rooted phylogenetic networks that were initially introduced by Cardona et al.~\cite{cardona09} and are formally defined in the next section. Specifically, given a tree-child network $\cN$, we show that each non-essential arc of  $\cN$ is a reticulation arc of a particular family of subgraphs of $\cN$. We call these subgraphs  {\it caterpillar ladders}. Although a caterpillar ladder can have arbitrary size, it always has exactly two non-essential arcs, one of which can be deleted. Intuitively, once one of the two non-essential arcs is deleted, the other arc becomes essential in the resulting network. Moreover these arcs can be identified in $\cN$ in polynomial time. This implies that we can preprocess a tree-child network and eliminate its non-essential arcs quickly before running any algorithm that searches through the exponentially-sized multiset of rooted phylogenetic trees that are embedded in $\cN$. Despite this encouraging result, we also show that, in general, deciding if a given arc of an arbitrary phylogenetic network is non-essential is $\Pi_2^P$-complete. With NP-complete and coNP-complete problems being placed on the first level of the polynomial-time hierarchy, $\Sigma_2^P$-complete and $\Pi_2^P$-complete problems are placed on the second level of this hie\-rarchy. Generally speaking, problems that are complete for the second level of the polynomial-time hierarchy are more difficult than problems that are complete for the first level. For more information about the polynomial-time hierarchy and its complexity classes, we refer the interested reader to~\cite{garey79,stockmeyer76}.

The remainder of the paper is organized as follows. Section~\ref{sec:preliminaries} provides definitions and terminology for phylogenetic networks and makes the notion of a non-essential arc precise. We also show that, in considering which arcs of a phylogenetic network are non-essential, it suffices to consider only reticulation arcs. Additionally, Section~\ref{sec:preliminaries} establishes two preliminary lemmas concerning the arcs in tree-child networks. In Section~\ref{sec:ladders}, we formally introduce caterpillar ladders and their embeddings in tree-child networks. These ladders and two of their variants are the main tool in characterizing the non-essential arcs in tree-child networks, the main theorem of the paper. This theorem is stated in Section~\ref{sec:ladders} and proved in Section~\ref{sec:the-proof}. In Section~\ref{sec:algorithm}, we show that the number of arcs that need to be deleted from a tree-child network $\cN$ so that the resulting network $\cN'$ has only essential arcs and embeds the same set of trees as $\cN$ is fixed, and thus does not depend on the ordering in which non-essential arcs are deleted from $\cN$. Subsequently, in the same section, we present an algorithm to identify a single non-essential arc in $\cN$ and show how this algorithm can be used to obtain $\cN'$ from $\cN$ in time that is cubic in the number of leaves of $\cN$. Lastly, in Section~\ref{sec:hardness}, we establish the hardness result mentioned at the end of the last paragraph.

\section{Preliminaries}\label{sec:preliminaries}

This section introduces notation and terminology. Throughout the paper, $X$ denotes a non-empty finite set. \\

\noindent {\bf Phylogenetic networks.} A {\em rooted binary phylogenetic network $\cN$ on $X$} is a rooted acyclic directed graph with no parallel arcs satisfying the following three properties:
\begin{enumerate}[(i)]
\item the (unique) root has in-degree zero and out-degree two;
\item a vertex of out-degree zero has in-degree one, and the set of vertices with out-degree zero is $X$; and
\item all other vertices either have in-degree one and out-degree two, or in-degree two and out-degree one.
\end{enumerate}
For technical reasons, if $|X|=1$, then we additionally allow $\cN$ to consist of the single vertex in $X$. We use $V(\cN)$ to denote the vertex set of $\cN$. The vertices of $\cN$ of out-degree zero are called {\em leaves}, and  $X$ is referred to as the {\em leaf set} of $\cN$. Furthermore, a vertex of in-degree one and out-degree two is a {\em tree vertex}, while a vertex of in-degree two and out-degree one is a {\em reticulation}. An arc directed into a reticulation is called a {\em reticulation arc}. All other arcs are called {\em tree arcs}. Moreover, a reticulation arc $(u, v)$ of $\cN$ is a {\em shortcut} if $\cN$ has a directed path from $u$ to $v$ avoiding $(u, v)$. A {\em rooted binary phylogenetic $X$-tree} is a rooted binary phylogenetic network on $X$ with no reticulations. To ease reading, we will refer to a rooted binary phylogenetic network and a rooted binary phylogenetic tree as a {\it phylogenetic network} and a {\it phylogenetic tree}, respectively, since all such networks and trees in this paper are rooted and binary. 

Now, let $\cN$ be a phylogenetic network, and let $P$ be a directed path in $\cN$. If $(v_1,v_2,\ldots,v_n)$ is the sequence of vertices that $P$ traverses, then we denote $P$ by $v_1,v_2,\ldots,v_n$. Furthermore, a directed path $P=v_1,v_2,\ldots,v_n$ in $\cN$ with $n\geq 1$ is a called a {\it tree path} of $\cN$ if each vertex $v_j$ with $j\in\{2,3,\ldots,n\}$ is a tree vertex or a leaf. 

\begin{figure}
\center
\scalebox{1.2}{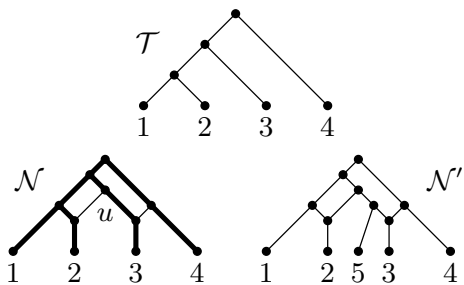}
\caption{Top: The caterpillar $\cT=(1,2,3,4)$. Bottom left: A phylogenetic network $\cN$ that is not tree-child and displays $\cT$. The thick arcs indicate an embedding of $\cT$ in $\cN$. Bottom right: A tree-child network $\cN'$.}
\label{fig:displaying}
\end{figure}

The main focus of this paper is on a particular class of phylogenetic networks which we now define. A phylogenetic network $\cN$ is called {\it tree-child} if each non-leaf vertex has a child that is a tree vertex or a leaf. Equivalently, $\cN$ is tree-child if, for each vertex $v$ in $\cN$ there exists a tree path from $v$ to a leaf. This equivalence is used frequently throughout this paper. To illustrate, Figure~\ref{fig:displaying} shows a non-tree-child network $\cN$ and a tree-child network $\cN'$ at the bottom left and bottom right, respectively. Note that $\cN$ is not tree-child, because both children of $u$ are reticulations.\\

\noindent {\bf Subtrees and caterpillars.} Let $\cT$ be a phylogenetic $X$-tree, and let $X'$ be a non-empty subset of $X$. Then $\cT|X'$ is the phylogenetic $X'$-tree obtained from the minimal  subtree of $\cT$ that connects all leaves in $X'$ by suppressing all vertices with in-degree one and out-degree one. We refer to $\cT|X'$ as the {\it restriction of $\cT$ to $X'$}. Furthermore, for two elements $\ell_1,\ell_2\in X$, we say that $\{\ell_1,\ell_2\}$ is a {\it cherry} of $\cT$ if $\ell_1$ and $\ell_2$ have the same parent.

Let $X=\{\ell_1,\ell_2,\ldots,\ell_n\}$ with $n\geq 2$, and let $\cT$ be a phylogenetic $X$-tree. For each $\ell_j\in X$ with $j\in\{1,2,\ldots,n\}$, let $p_j$ be the parent of $\ell_j$ in $\cT$.  We call $\cT$ a {\it caterpillar} if there exists an ordering, say $\ell_1,\ell_2,\ldots,\ell_n$, on the elements in $X$ such that, $p_1 = p_2$ and, for all $j \in \{2, 3,\ldots, n-1\}$ we have that $(p_{j+1}, p_j)$ is an arc in $\cT$. Furthermore, we denote such a caterpillar $\cT$ by $(\ell_1, \ell_2,\ell_3\ldots, \ell_n)$ or, equivalently, $(\ell_2, \ell_1, \ell_3, \ldots, \ell_n)$. An example of a caterpillar with cherry $\{1,2\}$ is shown in Figure~\ref{fig:displaying}.\\

\noindent {\bf Embeddings.} Let $\cN$ be a phylogenetic network on $X$, and let $\cT$ be a phylogenetic $X'$-tree, where $X'$ is a non-empty subset of $X$. Furthermore, let $\cE$ be a subset of the arcs of $\cN$.
We say that $\cN$ {\em displays} $\cT$ if $\cT$ can be obtained from $\cN$ by deleting arcs and vertices, and suppressing any resulting vertices of in-degree one and out-degree one. Moreover, $\cE$ is an {\it embedding of $\cT$} in $\cN$ if the subgraph of $\cN$ induced by the elements in $\cE$ is a subdivision of $\cT$.
If $\cE$ is an embedding of $\cT$ in $\cN$, then $\cT$ is displayed by $\cN$. Now, suppose that $\cE$ is an embedding of $\cT$ in $\cN$, and let $\{e_1,e_2,\ldots,e_m\}$ be a subset of the arcs in $\cN$. If $\{e_1,e_2,\ldots,e_m\}\subseteq \cE$, we say that  $\cE$ {\em uses} $\{e_1,e_2,\ldots,e_m\}$. To ease reading, if $\cE$ uses a singleton $\{e_1\}$, we write that $\cE$ uses $e_1$. The phylogenetic tree $\cT$ that is shown in Figure~\ref{fig:displaying} is displayed by each of the two phylogenetic networks $\cN$ and $\cN'$ that are shown in the same figure. An embedding of $\cT$ in $\cN$ is indicated by the thick arcs of $\cN$.\\

\noindent {\bf Essential arcs.} Let $\cN$ be a phylogenetic network on $X$, and let $e$ be an arc of $\cN$. If there exists a phylogenetic $X$-tree $\cT$ displayed by $\cN$ such that every embedding of $\cT$ in $\cN$ uses $e$, we say that $e$ is {\it essential}. Otherwise, if, for each phylogenetic $X$-tree $\cT$ displayed by $\cN$, there exists an embedding of $\cT$ that avoids $e$, we say that $e$ is {\it non-essential}. To decide which arcs of $\cN$ are non-essential, it follows from the next lemma that it suffices to consider only the reticulation arcs of $\cN$.

\begin{lemma}
Let $\cN$ be a phylogenetic network on $X$, and let $e$ be a tree arc of $\cN$. If $e$ is non-essential, then there is a reticulation arc of $\cN$ that is non-essential.
\label{non-essential}
\end{lemma}

\begin{proof}
Let $e=(t, u)$, and suppose that $e$ is non-essential. Let $P$ be a maximal length tree path starting at $u$. If $P$ ends at a leaf, then every embedding of a phylogenetic $X$-tree displayed by $\cN$ uses $e$, a contradiction. Thus $P$ ends at a tree vertex, $v$ say, whose two children are both reticulations. Let $w$ be one of these reticulations. Now, let $\cT$ be a phylogenetic $X$-tree that is displayed by $\cN$. Since $e$ is non-essential, it follows that there exists an embedding of $\cT$ in $\cN$ that avoids $e$ and, consequently, it also avoids $(v,w)$. Therefore, as $e$ is non-essential, $(v, w)$ is also non-essential. This completes the proof of the lemma.\qed
\end{proof}

It follows from Lemma~\ref{non-essential} that if every reticulation arc of a phylogenetic network $\cN$ is essential, then every tree arc of $\cN$ is essential. Thus, for the remainder of the paper, when considering non-essential arcs of a phylogenetic network, we restrict our attention to reticulation arcs.

We end this section with two lemmas and an observation that establish basic properties of essential and non-essential arcs in tree-child networks. Let $e=(u, w)$ be a tree arc of a phylogenetic network $\cN$. If $\cN$ is tree-child, then there exists a tree path starting at $u$, traversing $e$, and ending at a leaf in $\cN$. The next lemma now follows from
\cite[Theorem 1.1]{semple16}. This lemma is used freely throughout this paper.

\begin{lemma}\label{l:embedding}
Let $\cN$ be a tree-child network on $X$, and let $\cE$ be a subset of the arcs of $\cN$. Then $\cE$ is an embedding of a phylogenetic $X$-tree in $\cN$ if and only if $\cE$ contains every tree arc of $\cN$ and $\cE$ contains, for each reticulation $v$ in $\cN$ with parents $u$ and $u'$, exactly one of $(u,v)$ and $(u',v)$.
\end{lemma}

\noindent As an aside, it follows from Lemma~\ref{l:embedding} that every tree arc of a tree-child network is essential.

Let $\cN$ be a tree-child network, and let $e=(u,v)$ be a reticulation arc of $\cN$. We denote by $\cN\setminus\{e\}$ the network obtained from $\cN$ by either deleting $e$ and suppressing $u$ and $v$ if $u$ is not the root of $\cN$, or deleting both arcs incident with $u$ and suppressing $v$ if $u$ is the root of $\cN$. It follows from~\cite[Lemma 7(i)]{doecker21} that $\cN$ is tree-child. Moreover, with
the definition of an essential arc in mind, we have the following observation.

\begin{observation}\label{ob}
Let $e$ be a reticulation arc of a tree-child network $\cN$ on $X$. Then the following two statements are equivalent:
\begin{enumerate}[{\rm (a)}]
\item $e$ is essential, and
\item $\cN$ displays a phylogenetic $X$-tree that is not  displayed by $\cN\setminus\{e\}$.
\end{enumerate}
\end{observation}

\begin{lemma}\label{l:no-shortcut}
Let $v$ be a reticulation of a tree-child network $\cN$ on $X$. Let $e$ and $e'$ be the reticulation arcs that are directed into $v$. If neither $e$ nor $e'$ is a shortcut, then $e$ and $e'$ are both essential.
\end{lemma}

\begin{proof}
Let $e=(u,v)$ and let $e'=(u',v)$. Since $\cN$ is tree-child and neither $e$ nor $e'$ is a shortcut, there are three distinct elements $\ell_1$, $\ell_2$, and $\ell_3$ in $X$ such that there exists a tree path from $u$ to $\ell_1$, a tree path from $u'$ to $\ell_2$, and a tree path from $v$ to $\ell_3$. Moreover, by Lemma~\ref{l:embedding}, $\cN$ displays phylogenetic $X$-trees $\cT$ and $\cT'$ such that $\cT|\{\ell_1,\ell_2,\ell_3\}=(\ell_1,\ell_3,\ell_2)$ and $\cT'|\{\ell_1,\ell_2,\ell_3\}=(\ell_2,\ell_3,\ell_1)$. Furthermore, it is easily seen that $\cT$ is not displayed by $\cN\setminus\{e\}$ and $\cT'$ is not displayed by $\cN\setminus\{e'\}$. It now follows by Observation~\ref{ob} that $e$ and $e'$ are both essential.\qed
\end{proof}

\section{Caterpillar Ladders and Statement of the Main Result}\label{sec:ladders}

\begin{figure}[t]
\center
\scalebox{1}{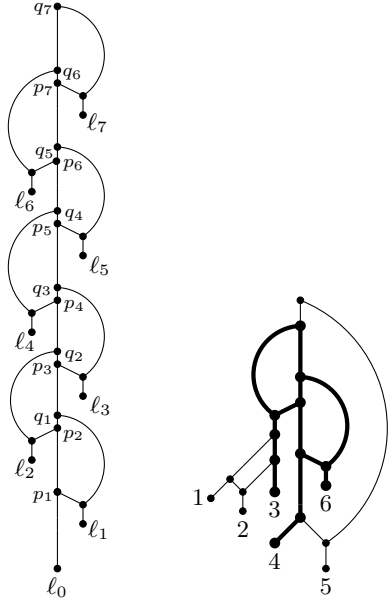}
\caption{Left: The caterpillar ladder $\langle\ell_0,\ell_1,\ell_2,\ldots,\ell_7\rangle$ with spine $q_7,q_6,p_7,q_5,p_6,q_4,p_5,q_3,p_4,q_2,p_3,q_1,p_2,p_1$. For each $j\in\{1,2,\ldots,7\}$, note that $v_j$ is the parent of $\ell_j$. Right: A tree-child network $\cN$ with thick arcs highlighting that $\langle4,6,3\rangle$ is a tight caterpillar ladder of $\cN$. Note that $\langle4,6,1\rangle$ is also a tight caterpillar ladder of $\cN$.}
\label{fig:caterpillar-ladder}
\end{figure}

Let $\cN$ be a phylogenetic network on $X=\{\ell_0,\ell_1,\ell_2,\ldots,\ell_k\}$, where $k\geq 1$, and let $\ell_0, \ell_1, \ell_2, \ldots, \ell_k$ be an ordering of the elements in $X$. Under this ordering, we say that $\cN$ is the (unique) {\it caterpillar ladder} on $X$ if, up to isomorphism, it has vertex set
$$X\bigcup_{j\in\{1,2,\ldots,k\}}\{v_j,p_j,q_j\},$$
and the vertices and arcs satisfy the following two properties:
\begin{enumerate}
\item [(i)] If $k=1$, then $(q_1,p_1)$ is an arc, where $q_1$ is the root of $\cN$. Otherwise, $q_k,q_{k-1},p_k,q_{k-2},p_{k-1},q_{k-3},p_{k-2},\ldots,q_3,p_4,q_2,p_3,q_1,p_2,p_1,\ell_0$ is a directed path, where $q_k$ is the root of $\cN$.

\item [(ii)] For each $j\in\{1,2,\ldots,k\}$, the vertex $v_j$ is a reticulation that is incident with the three arcs $(v_j,\ell_j)$, $(p_j,v_j)$, and $(q_j,v_j)$.
\end{enumerate}
For convenience, if $k=1$, we set $q_0=p_1$ and $p_2=q_1$. We denote the above caterpillar ladder by $\langle\ell_0,\ell_1,\ell_2,\ldots,\ell_k\rangle$. Furthermore, for $k>1$ (resp. $k=1$), we call the directed path $$q_k,q_{k-1},p_k,q_{k-2},p_{k-1},q_{k-3},p_{k-2},\ldots,q_3,p_4,q_2,p_3,q_1,p_2,p_1$$ (resp. the arc $(q_1,p_1))$ the {\it spine} of $\cN$. Note that a caterpillar ladder is tree-child. To illustrate, the caterpillar ladder $\langle\ell_0,\ell_1,\ell_2,\ldots,\ell_7\rangle$ is shown on the left-hand side of Figure~\ref{fig:caterpillar-ladder}.\\

\noindent {\bf Tight caterpillar ladders.} Let $\cN$ be a phylogenetic network on $X$, and let $\cC=\langle\ell_0,\ell_1,\ell_2,\ldots,\ell_k\rangle$ be a caterpillar ladder such that $\{\ell_0,\ell_1,\ell_2,\ldots,\ell_k\}\subseteq X$. We next define an embedding of $\cC$ in $\cN$. We say that $\cC$ is a {\it tight caterpillar ladder of $\cN$} if there exists a one-to-one map $\phi: V(\cC) \rightarrow V(\cN)$ with $\phi(\ell_j)=\ell_j$ for all $j\in \{0, 1, 2, \ldots, k\}$ satisfying the following three properties:\\

\noindent\begin{tabularx}{\linewidth}{@{}lX}
{\bf (P1)} &For $(p_1,\ell_0)$ and each arc $(v_j,\ell_j)$ in $\cC$ with $j\in\{1,2,\ldots,k\}$, there exists a tree path from $\phi(p_1)$ to $\phi(\ell_0)$ and from $\phi(v_j)$ to $\phi(\ell_j)$ in $\cN$, respectively.\\
{\bf (P2)} &For each reticulation arc $(u,v_j)$ in $\cC$ with $j\in\{1,2,\ldots,k\}$, there exists an arc $(\phi(u),\phi(v_j))$ in $\cN$. \\
{\bf (P3)}& For each arc $(u,w)$ on the spine of $\cC$, there exists an arc $(\phi(u),\phi(w))$ in $\cN$.
\end{tabularx}
\bigskip

To illustrate, the caterpillar ladder $\langle4,6,3\rangle$ is a tight caterpillar ladder of the phylogenetic network that is shown on the right-hand side of Figure~\ref{fig:caterpillar-ladder}. Now suppose that $\cC=\langle\ell_0,\ell_1,\ell_2,\ldots,\ell_k\rangle$ is a tight caterpillar ladder of $\cN$. Then it follows by (P3) and the existence of the two arcs $(\phi(q_k),\phi(q_{k-1}))$ and $(\phi(p_2),\phi(p_1))$ in $\cN$ that there exists no tight caterpillar ladder of $\cN$ whose leaf set strictly contains $\{\ell_0,\ell_1,\ell_2,\ldots,\ell_k\}$, that is, $\cC$ cannot be extended. Furthermore, in $\cN$, we refer to the arc $(\phi(p_1),\phi(v_1))$ as the {\it first rung} of $\cC$, and to the arc $(\phi(q_k),\phi(v_k))$ as the {\it last rung} of $\cC$.

We are now in a position to state the main result of this paper.

\begin{theorem}\label{t:main}
Let $e$ be a reticulation arc of a tree-child network $\cN$ on $X$. Then $e$ is non-essential if and only if $e$ is the first or last rung of a tight caterpillar ladder of $\cN$.
\end{theorem}

\begin{figure}
\center
\scalebox{1}{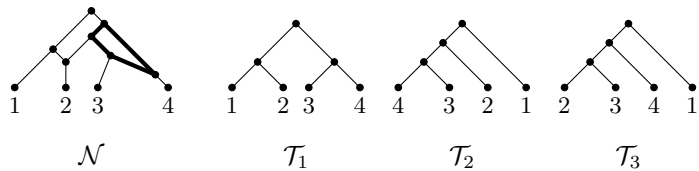}
\caption{A tree-child network $\cN$ with an avoidable cycle (indicated by thick arcs) that displays the three phylogenetic trees $\cT_1$, $\cT_2$, and $\cT_2$. }
\label{fig:avoidable-cycle}
\end{figure}

Before describing two embedding variants of caterpillar ladders in a phylogenetic network, we detail a connection between Theorem~\ref{t:main} and a previously established result. Cordue et al.~\cite{cordue14} have shown that a tree-child network $\cN$ on $X$ displays a phylogenetic $X$-tree twice precisely if $\cN$ contains an avoidable cycle. Roughly speaking, an underlying cycle $C$ of $\cN$ is {\it avoidable} if there exists an embedding $\cE$ of a phylogenetic $X$-tree $\cT$ in $\cN$ such that $\cE$ uses no more than two arcs $(u,v)$ of $\cN$ with the property that $u$ is a vertex of $C$ and $v$ is not a vertex of $C$, in which case $\cT$ is displayed twice by $\cN$.  Hence, if $\cN$ has a non-essential reticulation arc $e$, then $e$ is an arc of an avoidable cycle. Indeed, every underlying cycle of a tight caterpillar ladder in $\cN$ is avoidable. However, the converse does not hold. That is, if $\cN$ contains an avoidable cycle $C$, then $C$ does not necessarily contain a non-essential reticulation arc. Equivalently, if $\cN$ displays less than $2^k$ phylogenetic $X$-trees, where $k$ is the number of reticulations in $\cN$, then $\cN$ may or may not contain a non-essential reticulation arc. An example is shown in Figure~\ref{fig:avoidable-cycle}, where the phylogenetic tree $\cT_1$ is displayed twice by the tree-child network $\cN$ that is shown in the same figure. However, a straightforward check shows that all four reticulation arcs of $\cN$ are essential.

As an immediate consequence of Theorem~\ref{t:main}, we also have the following corollary. Let $\cN$ be a phylogenetic network. Then $\cN$ is called {\it normal}  if it is tree-child and has no shortcut. Moreover, $\cN$ is called  {\it level-$1$} if every underlying cycle of $\cN$ contains exactly one reticulation, that is, if the underlying cycles of $\cN$ are vertex disjoint. The classes of normal and level-$1$ networks are proper subclasses of tree-child networks. Note that if $\cN$ is normal, then it has no tight caterpillar ladder. Furthermore, if $\cC$ is a tight caterpillar ladder of a level-$1$ network, then $\cC$ has exactly one reticulation.

\begin{corollary}
Let $\cN$ be a phylogenetic network. If $\cN$ is normal, or if $\cN$ is level-$1$ and every underlying cycle has length at least four, then every reticulation arc of $\cN$ is essential.
\label{normal}
\end{corollary}

\noindent For normal networks, Corollary~\ref{normal} was independently established in~\cite{ier10} and~\cite{wil12} while, for level-$1$ networks, Corollary~\ref{normal} follows from~\cite{cordue14} as all of the underlying cycles of a level-$1$ network are unavoidable.

The remainder of this section contains  definitions of two variants of embeddings of caterpillar ladders in a phylogenetic network. These variants are used to establish Theorem~\ref{t:main} in the next section.\\

\noindent {\bf Nearly-tight and loose caterpillar ladders.} Let $\cN$ be a phylogenetic network on $X$, and let $\cC=\langle\ell_0,\ell_1,\ell_2,\ldots,\ell_k\rangle$ be a caterpillar ladder such that $\{\ell_0,\ell_1,\ell_2,\ldots,\ell_k\}\subseteq X$. If $k=1$, recall that $q_0=p_1$ and $p_2=q_1$. Furthermore, let $\phi: V(\cC) \rightarrow V(\cN)$ be a one-to-one map such that $\phi(\ell_j)=\ell_j$ for all $j\in \{0, 1, 2, \ldots, k\}$. First, we say that $\cC$ is a {\it nearly-tight caterpillar ladder of $\cN$} if $\phi$ satisfies (P1) and (P2), as well as exactly one of the following two weaker versions of (P3):\\

\noindent\begin{tabularx}{\linewidth}{@{}lX}
{\bf (P3\textsuperscript{$+$})} &For each arc $(u,w)$ on the spine of  $\cC$ with $(u,w)\ne (q_k,q_{k-1})$, there is an arc $(\phi(u),\phi(w))$ in $\cN$, and there exists a directed path from $\phi(q_k)$ to $\phi(q_{k-1})$ in $\cN$.\\
{\bf (P3\textsuperscript{$-$})}& For each arc $(u,w)$ on the spine of  $\cC$ with $(u,w)\ne (p_2,p_1)$, there is an arc $(\phi(u),\phi(w))$ in $\cN$, and there exists  a directed path from $\phi(p_2)$ to $\phi(p_1)$ in $\cN$.
\end{tabularx}
\bigskip

\noindent If (P3\textsuperscript{$+$}) applies, we refer to $\cC$ as $\cC_k^+$ and, if (P3\textsuperscript{$-$}) applies, we refer to $\cC$ as $\cC_k^-$. To clarify, if $k=1$ and $\cC$ is a nearly-tight caterpillar ladder of $\cN$, then there is a directed path from $q_1$ to $p_1$. Second, we say that $\cC$ is a {\it loose caterpillar ladder of $\cN$} if $\phi$ satisfies (P1) and (P2), as well as exactly one of the following two weaker versions of (P3):\\

\noindent\begin{tabularx}{\linewidth}{@{}lX}
{\bf (P3$^\uparrow$)} & For each arc $(u,w)$ on the spine of $\cC$ with $(u,w)\ne (q_k,q_{k-1})$, there exists a tree path from $\phi(u)$ to $\phi(w)$ in $\cN$, and there exists a directed path from $\phi(q_k)$ to $\phi(q_{k-1})$ in $\cN$.\\
{\bf (P3$^\downarrow$)} & For each arc $(u,w)$ on the spine of $\cC$ with $(u,w)\ne (p_2,p_1)$, there exists a tree path from $\phi(u)$ to $\phi(w)$ in $\cN$, and there exists a directed path from $\phi(p_2)$ to $\phi(p_1)$ in $\cN$.
\end{tabularx}
\bigskip

\noindent If (P3$^\uparrow$) applies, we refer to $\cC$ as $\cC_k^\uparrow$ and, if (P3$^\downarrow$) applies, we refer to $\cC$ as $\cC_k^\downarrow$. If $k=1$, then the definitions of a nearly-tight caterpillar ladder and a loose caterpillar ladder coincide. It is an immediate consequence of the above definitions that each tight caterpillar ladder of $\cN$ is also a nearly-tight caterpillar ladder of $\cN$, and each nearly-tight caterpillar ladder of $\cN$ is also a loose caterpillar ladder of $\cN$.

\section{Proof of Theorem~\ref{t:main}}\label{sec:the-proof} 

In this section, we establish Theorem~\ref{t:main}. Lemma~\ref{l:one-direction} establishes one  direction of the theorem, whereas the converse is an amalgamation of Lemmas~\ref{l:bottom-up} and~\ref{l:top-down}. The converse is obtained by repeatedly applying one of the latter two lemmas to extend a nearly-tight caterpillar ladder $\cC_1^+=\langle \ell_0,\ell_1\rangle$ (resp. $\cC_1^-=\langle \ell_0,\ell_1\rangle$) of a tree-child network $\cN$ to a tight caterpillar ladder $\cC=\langle \ell_0,\ell_1,\ell_2,\ldots,\ell_k\rangle$ of $\cN$ by moving towards (resp. away from) the root of $\cN$. Lemmas~\ref{l:embedding-up}--\ref{l:tree-path-down} and Corollary~\ref{2down} establish certain properties of loose caterpillar ladders that are frequently used in the proofs of Lemmas~\ref{l:bottom-up} and~\ref{l:top-down}. Lastly, as a reminder to the reader, throughout this section we freely use Lemma~\ref{l:embedding}.

Let $\cN$ be a phylogenetic network, and let $\cC=\langle\ell_0,\ell_1,\ell_2,\ldots,\ell_k\rangle$ be a caterpillar ladder that is a tight caterpillar ladder of $\cN$, where $k\ge 1$. Then there exists a map $\phi$ from the vertex set of $\cC$ to the vertex set of $\cN$ that satisfies (P1)--(P3). To ease reading, we use the following notation throughout this section, where $j\in\{1,2,\ldots,k\}$:
\begin{enumerate}[1.]
\item For each reticulation $v_j$ in $\cC$, the vertex $\phi(v_j)$ in $\cN$ is denoted by $v_j$.
\item For each vertex $u$ on the spine of $\cC$, the vertex $\phi(u)$ in $\cN$ is denoted by $u$. Hence, each reticulation $v_j$ in $\cN$ has parents $p_j$ and $q_j$.
\item For each  $j$, the arc $(p_j,v_j)$ in $\cN$ is denoted by $e_j$.
\item For each  $j$, the arc $(q_j,v_j)$ in $\cN$ is denoted by $f_j$.
\end{enumerate}
Note that each $f_j$ is a shortcut of $\cN$. The notation naturally carries over for a nearly-tight caterpillar ladder and a loose caterpillar ladder of $\cN$.

The next lemma establishes one direction of Theorem~\ref{t:main}.

\begin{lemma}\label{l:one-direction}
Let $\cC=\langle\ell_0,\ell_1,\ell_2,\ldots,\ell_k\rangle$ be a tight caterpillar ladder of a tree-child network $\cN$ on $X$. Then each element in $\{e_1,f_k\}$ is non-essential and each element in $\{f_1,e_2,f_2,e_3,f_3,\ldots,e_{k-1},f_{k-1},e_k\}$ is essential. 
\end{lemma}

\begin{proof}
We first show that each of $e_1$ and $f_k$ is non-essential. Let $\cE$ be an embedding of a phylogenetic $X$-tree $\cT$ in $\cN$ that uses $e_1$, and let $j$ be the maximum element in $\{1,2,\ldots,k\}$ such that $\cE$ also uses $\{e_2,e_3,\ldots,e_j\}$. Observe that $\cT|\{\ell_0,\ell_1,\ell_2,\ldots,\ell_j\}=(\ell_0,\ell_1,\ell_2,\ldots,\ell_j)$. Now, as the choice of $j$ implies that $\cE$ uses $f_{j+1}$ unless $j=k$, it is straightforward to check that $$(\cE-\{e_1,e_2,\ldots,e_j\})\cup\{f_1,f_2,\ldots,f_j\}$$ is an embedding of $\cT$ in $\cN$. If $j=k$, then this set of arcs is also an embedding of $\cT$ in $\cN$. Thus $e_1$ is non-essential. Similarly, let $\cE'$ be an embedding of a phylogenetic $X$-tree $\cT'$ in $\cN$ that uses $f_k$, and let $j'$ be the minimum element in $\{1,2,\ldots,k\}$ such that $\cE'$ also uses $\{f_{j'},f_{j'+1},\ldots,f_{k-1}\}$. Then $\cT'|\{\ell_0,\ell_{j'},\ell_{j'+1},\ldots,\ell_k\}=(\ell_0,\ell_{j'},\ell_{j'+1},\ldots,\ell_k)$. Again, as the choice of $j'$ implies that $\cE'$ uses $e_{j'-1}$ unless $j'=1$, it follows that  $$(\cE'-\{f_{j'},f_{j'+1},\ldots,f_k\})\cup\{e_{j'},e_{j'+1},\ldots,e_k\}$$ is an embedding of $\cT'$ in $\cN$. If $j'=1$, then this set of arcs is also an embedding of $\cT$ in $\cN$. Thus $f_k$ is non-essential.
 
We complete the proof by showing that each element in $$\{f_1,e_2,f_2,e_3,f_3,\ldots,e_{k-1},f_{k-1},e_k\}$$ is essential, where $k\ge 2$. Let $e_j\in\{e_2,e_3,\ldots,e_k\}$, and let $\cE$ be an embedding of a phylogenetic $X$-tree $\cT$ in $\cN$ that uses $\{f_{j-1},e_j,\}$. Then $\cT|\{\ell_0,\ell_{j-1},\ell_j\}=(\ell_0,\ell_j,\ell_{j-1})$. However, for every phylogenetic $X$-tree $\cT'$ displayed by $\cN\setminus\{e_j\}$, we have $\cT'|\{\ell_0,\ell_{j-1},\ell_j\}=(\ell_0,\ell_{j-1},\ell_j)$. Hence, by Observation~\ref{ob}, $e_j$ is essential. Now let $f_j\in\{f_1,f_2,\ldots,f_{k-1}\}$, and let $\cE$ be an embedding of a phylogenetic $X$-tree $\cT$ in $\cN$ that uses $\{f_j,e_{j+1},\}$. Then $\cT|\{\ell_0,\ell_j,\ell_{j+1}\}=(\ell_0,\ell_{j+1},\ell_j)$. But, for every phylogenetic $X$-tree $\cT'$ displayed by $\cN\setminus\{f_j\}$, we have $\cT'|\{\ell_0,\ell_{j},\ell_{j+1}\}=(\ell_0,\ell_j,\ell_{j+1})$. Hence, again by Observation~\ref{ob}, $f_j$ is essential.\qed
\end{proof}

The next four lemmas and subsequent corollary establish properties of a loose caterpillar ladder of a tree-child network. These properties are then used in the proofs of  Lemmas~\ref{l:bottom-up} and~\ref{l:top-down} that, taken together, establish the converse of Theorem~\ref{t:main}.

\begin{lemma}\label{l:embedding-up}
Let $\cC_k^\uparrow=\langle\ell_0,\ell_1,\ell_2,\ldots,\ell_k\rangle$ be a loose caterpillar ladder of a tree-child network $\cN$ on $X$.
If $\cE$ and $\cE'$ are two embeddings of the same phylogenetic $X$-tree $\cT$ in $\cN$ such that $\cE$ uses $\{e_1,e_2,\ldots,e_k\}$ and $\cE'$ uses $f_1$, then $\cE'$ also uses $\{f_2, f_3, \ldots, f_k\}$. 
\end{lemma}

\begin{proof}
If $k=1$, then the lemma vacuously holds, so we may assume that $k\ge 2$. By  (P3$^{\uparrow}$), there is a unique tree path from $q_{k-1}$ to $p_1$ in $\cN$. Hence, since $\cE$ uses $\{e_1,e_2,\ldots,e_k\}$, we have $\cT|\{\ell_0,\ell_1,\ell_2,\ldots,\ell_k\}=(\ell_0,\ell_1,\ell_2,\ldots,\ell_k)$. Now, towards a contradiction, assume that $\cE'$ uses $e_j$ for some $j\in\{2,3,\ldots,k\}$. Without loss of generality, we may assume that $j$ is maximized so that $\cE'$ uses $\{f_1,f_2,\ldots,f_{j-1}\}$. Then, by considering $\cE'$, we have $\cT|\{\ell_0,\ell_1,\ell_2,\ldots,\ell_j\}=(\ell_0,\ell_1,\ell_2,\ldots,\ell_{j-2},\ell_j,\ell_{j-1})$ which gives a contradiction since $\cE$ and $\cE'$ are embeddings of the same phylogenetic $X$-tree, and $j\ge 2$.\qed
\end{proof}

Note that the analogue of the next lemma for when $k=2$ is established as Corollary~\ref{2down}.

\begin{lemma}\label{l:embedding-down}
Let $\cC_k^\downarrow=\langle\ell_0,\ell_1,\ell_{2},\ldots,\ell_k\rangle$ be a loose caterpillar ladder of a tree-child network $\cN$ on $X$, where $k=1$ or $k\ge 3$. If $\cE$ and $\cE'$ are two embeddings of the same phylogenetic $X$-tree $\cT$ in $\cN$ such that $\cE'$ uses $\{f_1,f_{2},\ldots,f_k\}$ and $\cE$ uses $e_k$, then $\cE$ also uses $\{e_1, e_{2}, \ldots, e_{k-1}\}$.
\end{lemma}

\begin{proof}
If $k=1$, then the lemma vacuously holds. Thus assume that $k\ge 3$. By (P3$^{\downarrow}$), there is a unique tree path from $q_k$ to $p_2$ in $\cN$. Therefore, as $\cE'$ uses $\{f_1,f_{2},\ldots,f_k\}$, it follows that  $\cT|\{\ell_1,\ell_{2},\ldots,\ell_k\}=(\ell_1,\ell_{2},\ldots,\ell_k)$. Now, towards a contradiction, assume that $\cE$ uses $f_j$ for some $j\in\{1,2,\ldots,k-1\}$. Without loss of generality, we may assume that $j$ is minimized so that $\cE$ uses $\{e_{j+1},e_{j+2},\ldots,e_k\}$. If $j\leq k-2$, then, by considering $\cE$, it follows that $\cT|\{\ell_j,\ell_{j+1},\ell_{j+2}\ldots,\ell_k\}=(\ell_{j+1},\ell_j,\ell_{j+2},\ldots,\ell_k)$, a contradiction as $j\le k-2$. Otherwise, if $j=k-1$, then, regardless of whether $\cE$ uses $e_{k-2}$ or $f_{k-2}$, we have $\cT|\{\ell_{k-2},\ell_{k-1},\ell_{k}\}=(\ell_{k-2},\ell_k,\ell_{k-1})$, again a contradiction as $k\ge 3$.\qed
\end{proof}

\begin{lemma}\label{l:tree-path-up} 
Let $\cC_k^\uparrow=\langle\ell_0,\ell_1,\ell_2,\ldots,\ell_k\rangle$ be a loose caterpillar ladder of a tree-child network $\cN$ on $X$. Suppose that $e_1$ is a non-essential arc in $\cN$.~Then
\begin{enumerate}[{\rm (i)}]
\item there is a unique tree path from $q_k$ to $q_{k-1}$ in $\cN$, and

\item every tree path in $\cN$ that starts at $q_k$ traverses $p_1$.
\end{enumerate}
\end{lemma}

\begin{proof}
Throughout the proof, let $\cE$ be an embedding of a phylogenetic \mbox{$X$-tree} $\cT$ in $\cN$ that uses $\{e_1,e_2,\ldots,e_k\}$. We first show that there is a tree path from $q_k$ to $q_{k-1}$ from which the uniqueness result is immediate. Let $P$ be a directed path from $q_k$ to $q_{k-1}$. By (P3$^\uparrow$), it follows that $P$ exists. Assume that $P$ is not a tree path. Let $w$ be the first reticulation on $P$, and let $u$ be the parent of $w$ that lies on $P$. Since $\cN$ is tree-child and $q_k$ is a parent of a reticulation, it follows that $u$ is a tree vertex and $u\ne q_k$. Moreover, there exists a tree path from $u$ to a leaf $x$ with $x\notin\{\ell_0,\ell_1,\ldots,\ell_k\}$, and so
\begin{equation}\label{eq:tree1}
\cT|\{x,\ell_0,\ell_1,\ldots,\ell_k\}=(\ell_0,\ell_1,\ldots,\ell_k,x). 
\end{equation}
Since $e_1$ is non-essential, there exists an embedding $\cE'$ of $\cT$ in $\cN$ such that $\cE'$ uses $f_1$. By Lemma~\ref{l:embedding-up}, $\cE'$ also uses $\{f_2,f_3,\ldots,f_k\}$. If $\cE'$ uses $(u,w)$, then $\cT|\{x,\ell_0,\ell_1,\ldots,\ell_k\}=(\ell_0,\ell_1,\ldots,\ell_{k-1},x,\ell_k)$. On the other hand, if $\cE'$ does not $(u,w)$, it  follows that $\cT|\{x,\ell_0,\ell_1,\ldots,\ell_k\}$ is $(\ell_0,\ell_1,\ldots,\ell_{k-1},x,\ell_k)$ or it has $\{\ell_k,x\}$ as a cherry. All outcomes contradict the restriction of $\cT$ in~(\ref{eq:tree1}). Hence $P$ is the unique tree path from $q_k$ to $q_{k-1}$. 

We complete the proof by showing that every tree path that starts at $q_k$ traverses $p_1$. Assume that there is a tree path $Q$ that starts at $q_k$ and does not traverse $p_1$. Let $P'$ be the (unique) tree path in $\cN$ that starts at $q_k$, ends at $\ell_0$, and traverses in order
$$q_{k-1},p_k,q_{k-2},p_{k-1},\ldots,q_2,p_3,q_1,p_2,p_1.$$
Furthermore, let $w$ be the first vertex on $Q$ that does not lie on $P'$. As $\cN$ is tree-child, there exists a tree path from $w$ to a leaf $x$ such that $x\notin\{\ell_0,\ell_1,\ldots,\ell_k\}$. Now consider the maximum length subpath of $Q$ that coincides with $P'$, and let $u$ be the last vertex of that subpath that is an element in $\{q_1,p_2,q_2,p_3,q_3,\ldots,p_k,q_k\}$. Since $Q$ does not traverse $p_1$, we have $u\ne p_1$. We next distinguish two cases. First, assume that $u=p_j$ for some $j\in\{2,3,\ldots,k\}$. Then, $k\ge 2$ and, by considering $\cE$,
$$\cT|\{x,\ell_0,\ell_1,\ldots,\ell_k\}=(\ell_0,\ell_1, \ldots,\ell_{j-1},x,\ell_j,\ldots,\ell_k).$$
Since $e_1$ is non-essential, there exists an embedding $\cE'$ of $\cT$ in $\cN$ such that $\cE'$ uses $f_1$ and, by Lemma~\ref{l:embedding-up}, it also uses $\{f_2,f_3,\ldots,f_k\}$. Considering $\cE'$, it now follows that $$\cT|\{x,\ell_0,\ell_1,\ldots,\ell_k\}=(\ell_0,\ell_1,\ldots,x,\ell_{j-1},\ell_j,\ldots,\ell_k),$$ a contradiction. Thus $u\neq p_j$. Second, assume that $u=q_j$ for some $j\in\{1,2,\ldots,k\}$. Then, by considering $\cE$, 
$$\cT|\{x,\ell_0,\ell_1,\ldots,\ell_k\}=(\ell_0,\ell_1,\ldots,\ell_j,\ell_{j+1},x,\ldots,\ell_k)$$
if $u\ne q_k$ or
$$\cT|\{x,\ell_0,\ell_1,\ldots,\ell_k\}=(\ell_0,\ell_1,\ldots,\ell_k,x)$$
if $u=q_k$.
Again, since $e_1$ is non-essential, there exists an embedding $\cE'$ of $\cT$ in $\cN$ such that $\cE'$ uses $f_1$. By Lemma~\ref{l:embedding-up}, $\cE'$ also uses $\{f_2,f_3,\ldots,f_k\}$. Regardless of whether or not $u=q_k$, it now follows that
$$\cT|\{x,\ell_0,\ell_1,\ldots,\ell_k\}=(\ell_0,\ell_1,\ldots,\ell_{j-1},x,\ell_j,\ldots,\ell_k),$$ 
another contradiction. This establishes the lemma. \qed
\end{proof}

For the next lemma, recall that if $k=1$, we set $p_2=q_1$.

\begin{lemma}\label{l:tree-path-down}
Let $\cC_k^\downarrow=\langle\ell_0,\ell_1,\ell_{2},\ldots,\ell_k\rangle$ be a loose caterpillar ladder of a  tree-child network $\cN$ on $X$. Suppose that $f_k$ is a non-essential arc in $\cN$.~Then
\begin{enumerate}[{\rm (i)}]
\item there is a unique tree path from $p_2$ to $p_1$ in $\cN$, and

\item every tree path in $\cN$ that starts at $p_2$ traverses $p_1$.
\end{enumerate}
\end{lemma}

\begin{proof}
The proof is similar to that of Lemma~\ref{l:tree-path-up}. Throughout the proof, let $\cE'$ be an embedding of a phylogenetic $X$-tree $\cT$ in $\cN$ that uses $\{f_1,f_2,\ldots,f_k\}$. 

To prove (i), let $P$ be a directed path from $p_2$ to $p_1$. By (P3$^{\downarrow}$), such a path exists. Assume that $P$ is not a tree path. Let $w$ be the first reticulation on $P$, and let $u$ be the parent of $w$ that lies on $P$. Since $\cN$ is tree-child and $p_2$ is a parent of a reticulation, it follows that $u$ is a tree vertex and $u\ne p_2$. Moreover, there exists a tree path from $u$ to a leaf $x$ with $x\notin\{\ell_0,\ell_1,\ell_{2}\ldots,\ell_k\}$.  By considering $\cE'$, it is easily checked that $\cT|\{x, \ell_0, \ell_1\}$ does not have $\{\ell_0, \ell_1\}$ as a cherry. Since $f_k$ is non-essential, there exists an embedding $\cE$ of $\cT$ in $\cN$ such that $\cE$ uses $e_k$. If $k=1$ or $k\ge 3$, then, by Lemma~\ref{l:embedding-down}, $\cE$ also uses $\{e_1,e_{2},\ldots,e_{k-1}\}$. Regardless of whether $\cE$ uses $(u,w)$ or not, it  follows that $\cT|\{x,\ell_0,\ell_1\}$ has $\{\ell_0,\ell_1\}$ as a cherry, a contradiction. So assume that $k=2$. Recall that $\cE'$ uses $f_2$ and $\cE$ uses $e_2$. If $\cE$ uses $e_1$, then $\cT|\{x, \ell_0, \ell_1, \ell_2\}$ has $\{\ell_0, \ell_1\}$ as a cherry, a contradiction. If $\cE$ uses $f_1$, then $\cT|\{x, \ell_1, \ell_2\}=(x, \ell_2, \ell_1)$ but, by considering $\cE'$, we have $\cT|\{x, \ell_1, \ell_2\}=(x, \ell_1, \ell_2)$, another contradiction. Hence $P$ is the unique tree path from $p_2$ to $p_1$, and so (i) holds.

To complete the proof, assume that there exists a tree path $Q$ that starts at $p_{2}$ and does not traverse $p_1$. We use an argument that is similar to that of the first part of the proof.  Let $w$ be the first vertex on $Q$ that does not lie on the unique tree path from $p_2$ to $p_1$ in $\cN$. Then there exists a tree path from $w$ to a leaf $x$ such that $x\notin\{\ell_0,\ell_1,\ell_{2},\ldots,\ell_k\}$. Since (i) holds, it follows by considering $\cE'$ that
\begin{equation}\label{eq:treeX}
\cT|\{x, \ell_0, \ell_1\}=(\ell_0, x, \ell_1).
\end{equation}
Since $f_k$ is non-essential, there exists an embedding $\cE$ of $\cT$ in $\cN$ such that $\cE$ uses $e_k$. If $k=1$ or $k\ge 3$, then, by Lemma~\ref{l:embedding-down}, it also uses $\{e_1,e_{2},\ldots,e_{k-1}\}$. Thus
$$\cT|\{x, \ell_0, \ell_1\}=(\ell_0, \ell_1, x),$$
contradicting (\ref{eq:treeX}).
If $k=2$ and $\cE$ uses $e_1$, then $\cT|\{x, \ell_0, \ell_1, \ell_2\}=(\ell_0, \ell_1, x, \ell_2)$, a contradiction to the restriction of $\cT$ in~(\ref{eq:treeX}). Furthermore, if $k=2$ and $\cE$ uses $f_1$, then $\cT|\{x, \ell_0, \ell_1, \ell_2\}=(\ell_0, x, \ell_2, \ell_1)$, but, by considering $\cE'$, we have $\cT|\{x,\ell_0, \ell_1, \ell_2\}=(\ell_0, x,\ell_1,\ell_2)$. This last contradiction completes the proof of the lemma.\qed
\end{proof}

The next corollary settles Lemma~\ref{l:embedding-down} for when $k=2$ with an additional assumption that $f_2$ is non-essential.

\begin{corollary}
Let $\cC^{\downarrow}_2=\langle\ell_0, \ell_1, \ell_2\rangle$ be a loose caterpillar ladder of a tree-child network $\cN$ on $X$.
Suppose that $f_2$ is non-essential. If $\cE$ and $\cE'$ are two embeddings of the same phylogenetic tree $\cT$ in $\cN$ such that $\cE'$ uses $\{f_1, f_2\}$ and $\cE$ uses $e_2$, then $\cE$ also uses $e_1$.
\label{2down}
\end{corollary}

\begin{proof}
It follows by Lemma~\ref{l:tree-path-down} that the directed path from $p_2$ to $p_1$ in the definition of $\cC^{\downarrow}_2$ is a tree path. Hence, by considering $\cE'$, we have $\cT|\{\ell_0, \ell_1, \ell_2\}=(\ell_0, \ell_1, \ell_2)$. Thus if $\cE$ does not use $e_1$, then, by considering $\cE$, we have
$$\cT|\{\ell_0, \ell_1, \ell_2\}=(\ell_0, \ell_2, \ell_1),$$
a contradiction. Hence $\cE$ uses $\{e_1, e_2\}$.\qed
\end{proof}

\begin{lemma}\label{l:bottom-up}
Let $\cC_k^+=\langle\ell_0,\ell_1,\ell_2,\ldots,\ell_k\rangle$ be a nearly-tight caterpillar ladder of a tree-child network $\cN$ on $X$.
If $e_1$ is non-essential in $\cN$, then exactly one of the following two properties holds:
\begin{enumerate}[{\rm (i)}]
\item $\cC_k^+$ is a tight caterpillar ladder of $\cN$, or
\item there exists an element $\ell_{k+1}\in X-\{\ell_0,\ell_1,\ell_2,\ldots,\ell_k\}$ such that $\cC_{k+1}^+=\langle\ell_0,\ell_1,\ell_2,\ldots,\ell_k,\ell_{k+1}\rangle$ is a nearly-tight caterpillar ladder of $\cN$. 
\end{enumerate} 
\end{lemma}

\begin{proof}
Assume that (i) does not hold, that is $(q_k,q_{k-1})$ is not an arc in $\cN$. We will show that (ii) holds. By Lemma~\ref{l:tree-path-up}, it follows that there is a unique tree path $P$ from $q_k$ to $q_{k-1}$. Now, let $p_{k+1}$ be the child of $q_k$ such that $p_{k+1}\ne v_k$. Since $P$ is a tree path, $p_{k+1}$ is a tree vertex and $p_{k+1}$ has a child $v_{k+1}$ that does not lie on $P$.  Applying  Lemma~\ref{l:tree-path-up} again, it follows that every tree path that starts at $q_k$ traverses $p_1$. Hence, as there is a unique tree path from $q_k$ to $p_1$, it follows that $v_{k+1}$ is a reticulation. Let $q_{k+1}$ be the parent of $v_{k+1}$ such that $q_{k+1}\ne p_{k+1}$. Furthermore, let $\ell_{k+1}$  be a leaf at the end of a tree path that starts at $v_{k+1}$. As $\cN$ is tree-child, $\ell_{k+1}\notin \{\ell_0,\ell_1,\ldots,\ell_k\}$.

Following the notation of a nearly-tight caterpillar ladder, let $e_{k+1}=(p_{k+1},v_{k+1})$, and let $f_{k+1}=(q_{k+1},v_{k+1})$. Furthermore, let $\cE$ be an embedding of a phylogenetic $X$-tree $\cT$ in $\cN$ that uses $\{e_1,e_2,\ldots,e_k,e_{k+1}\}$. Since $\cT|\{\ell_0,\ell_1,\ldots,\ell_k\}=(\ell_0,\ell_1,\ldots,\ell_k)$, it follows that
\begin{equation}\label{eq:leading-tree}
\cT|\{\ell_0,\ell_1,\ldots,\ell_k,\ell_{k+1}\}=(\ell_0,\ell_1,\ldots,\ell_k,\ell_{k+1}).
\end{equation}
As $e_1$ is non-essential, there exists an embedding $\cE'$ of $\cT$ in $\cN$ that uses $f_1$. By Lemma~\ref{l:embedding-up}, $\cE'$ also uses $\{f_2,f_3,\ldots,f_k\}$ and, by the restriction of $\cT$ in (\ref{eq:leading-tree}), $\cE'$ uses $f_{k+1}$.

\begin{sublemma}\label{s:shortcut-up}
The arc $f_{k+1}$ is a shortcut.  
\end{sublemma}

\begin{proof}
First assume that neither $e_{k+1}$ nor  $f_{k+1}$ is a shortcut, in which case $q_{k+1}$ does not lie on $P$. Let $x$  be a leaf at the end of a tree path that starts at $q_{k+1}$. As $\cN$ is tree-child, $x\ne\ell_{k+1}$. Furthermore, as $e_{k+1}$ and $f_{k+1}$ are not shortcuts, $x\notin\{\ell_0,\ell_1,\ldots,\ell_k\}$. Therefore, by considering $\cE$, it follows that $\cT|\{x,\ell_0,\ell_1,\ldots,\ell_k,\ell_{k+1}\}=(\ell_0,\ell_1,\ldots,\ell_k,\ell_{k+1},x)$. But, by considering $\cE'$ and recalling that $\cE'$ uses $f_{k+1}$, it is easily checked that $\cT|\{x,\ell_0,\ell_1,\ldots,\ell_k,\ell_{k+1}\}$ has $\{x,\ell_{k+1}\}$ as a cherry, a contradiction. Hence, $x\in\{\ell_0,\ell_1,\ldots,\ell_k\}$, and so one of $e_{k+1}$ and $f_{k+1}$ is a shortcut. Assume that $e_{k+1}$ is a shortcut. Then, by considering $\cE$, we again have the restriction of $\cT$ in (\ref{eq:leading-tree}). On the other hand, as $\cC_k^+$ is nearly-tight, by considering $\cE'$, one of the following statements holds:
\begin{enumerate}[(I)]
\item If $x=\ell_0$, then $q_{k+1}$ lies on the subpath of $P$ from $p_{k+1}$ to $q_{k-1}$ and $\cT|\{\ell_0,\ell_1,\ldots,\ell_k,\ell_{k+1}\}=(\ell_0,\ell_1,\ldots,\ell_{k-1},\ell_{k+1},\ell_k)$, or $q_{k+1}$ lies on the tree path from $p_1$ to $\ell_0$ and $\cT|\{\ell_0,\ell_1,\ldots,\ell_k,\ell_{k+1}\}$ has $\{\ell_0,\ell_{k+1}\}$ as a cherry.
\item If $x\in\{\ell_1,\ell_2,\ldots,\ell_k\}$, then $q_{k+1}$ lies on a tree path from $v_j$ to $\ell_j$ for some $j\in\{1,2,\ldots,k\}$ and $\cT|\{\ell_0,\ell_1,\ldots,\ell_k,\ell_{k+1}\}$ has $\{\ell_j,\ell_{k+1}\}$ as a cherry.
\end{enumerate}
Both cases contradict~(\ref{eq:leading-tree}); thereby implying that $f_{k+1}$ is a shortcut.\qed
\end{proof}

By the existence of $P$ and as $\cC_k^+$ is a nearly-tight caterpillar ladder of $\cN$, we have that
\begin{sublemma}\label{loose-up}
$\cC_{k+1}^\uparrow=\langle\ell_0,\ell_1,\ell_2,\ldots,\ell_k,\ell_{k+1}\rangle$ is a loose caterpillar ladder of $\cN$.
\end{sublemma}
\noindent Now, applying Lemma~\ref{l:tree-path-up} to $\cC_{k+1}^\uparrow$ implies that there is a unique tree path~$P'$ from $q_{k+1}$ to $q_k$. 

\begin{sublemma}
Let $u$ be the child of $p_{k+1}$ on $P$. Then $u=q_{k-1}$.
\label{child1}
\end{sublemma}

\begin{proof}
Assume that $u\ne q_{k-1}$. By Lemma~\ref{l:tree-path-up}, every tree path that starts at $q_{k+1}$ traverses $p_1$. Hence, the child $w$ of $u$ that does not lie on $P$ is a reticulation. Let $u'$ be the parent of $w$ that is not $u$, and let $x$ be a leaf at the end of a tree path that starts at $w$. Since $\cN$ is tree-child, $x\notin\{\ell_0,\ell_1,\ell_2,\ldots,\ell_{k+1}\}$. As before, we  play several embeddings of phylogenetic trees off against each other to show that $u=q_{k-1}$. Specifically, let $\cE_1$ be an embedding of a phylogenetic $X$-tree $\cT_1$ in $\cN$ that uses $\{(u,w),e_1,e_2,\ldots,e_k,e_{k+1}\}$. Then 
\begin{equation}\label{eq:leading-tree-two}
\cT_1|\{x,\ell_0,\ell_1,\ldots,\ell_k,\ell_{k+1}\}=(\ell_0,\ell_1,\ldots,\ell_k,x,\ell_{k+1}).
\end{equation}
As $e_1$ is non-essential, there exists an embedding $\cE_1'$ of $\cT_1$ in $\cN$ that uses $f_1$. Furthermore, by Lemma~\ref{l:embedding-up} and (\ref{loose-up}), $\cE'_1$ also uses $\{f_2,f_3,\ldots,f_k,f_{k+1}\}$. Lastly, by the restriction of $\cT_1$ in (\ref{eq:leading-tree-two}), $\cE_1'$ uses $(u',w)$
and there is a directed path $Q$ from a vertex $s$ on $P'$ to $w$ that avoids $p_{k+1}$ and traverses $(u',w)$. Assume that $s\ne u'$. Let $y$ be the leaf at the end of tree path that starts at $u'$. As $\cN$ is tree-child, $y\ne x$. Furthermore, since $Q$ avoids $p_{k+1}$, and $P$ is a tree path, it is easily seen that $y\not\in \{\ell_0,\ell_1,\ldots,\ell_{k+1}\}$. Turning back to $\cT_1$, it follows that $\cT_1|\{x,y,\ell_0,\ell_1,\ldots,\ell_k,\ell_{k+1}\}$ does not have $\{x,y\}$ as a cherry by considering $\cE_1$ and does have $\{x,y\}$ as a cherry by considering $\cE_1'$, a contradiction. Hence $s=u'$ and so $(u',w)$ is a shortcut.

Now let $\cE_2$ be an embedding of a phylogenetic $X$-tree $\cT_2$ in $\cN$ that uses $\{(u',w),e_1,e_2,\ldots,e_k,e_{k+1}\}$. Then 
\begin{equation}\label{eq:leading-tree-three}
\cT_2|\{x,\ell_0,\ell_1,\ldots,\ell_k,\ell_{k+1}\}=(\ell_0,\ell_1,\ldots,\ell_{k-1},\ell_k,\ell_{k+1},x). 
\end{equation}
Since $e_1$ is non-essential, there exists an embedding $\cE_2'$ of $\cT_2$ in $\cN$ that uses $f_1$. By Lemma~\ref{l:embedding-up} and (\ref{loose-up}), $\cE_2'$ also uses~$\{f_2,f_3,\ldots,f_k,f_{k+1}\}$. 
We continue to consider $\cE_2'$. If $\cE_2'$ uses $(u,w)$, then $$\cT_2|\{x,\ell_0,\ell_1,\ldots,\ell_k,\ell_{k+1}\}=(\ell_0,\ell_1,\ldots,\ell_{k-1},x,\ell_k,\ell_{k+1})$$ and, if  $\cE_2'$ uses $(u',w)$, then $$\cT_2|\{x,\ell_0,\ell_1,\ldots,\ell_k,\ell_{k+1}\}=(\ell_0,\ell_1,\ldots,\ell_{k-1},\ell_k,x,\ell_{k+1}).$$ Both cases contradict that $\cE_2$ and $\cE_2'$ are embeddings of $\cT_2$. Hence $u=q_{k-1}$.\qed
\end{proof}

It now follows by (\ref{loose-up}) and (\ref{child1}) that $\cC_{k+1}^+=\langle\ell_0,\ell_1,\ell_2,\ldots,\ell_k,\ell_{k+1}\rangle$ is a nearly-tight caterpillar ladder of $\cN$. This completes the proof of the lemma.\qed
\end{proof}

\begin{lemma}\label{l:top-down}
Let $\cC_{k-1}^-=\langle\ell_0,\ell_2,\ell_{3},\ldots,\ell_k\rangle$ be a nearly-tight caterpillar ladder of a tree-child network $\cN$ on $X$.
If $f_k$ is non-essential in $\cN$, then exactly one of the following two properties holds:
\begin{enumerate}[{\rm (i)}]
\item $\cC_{k-1}^-$ is a tight caterpillar ladder of $\cN$, or
\item there exists an element $\ell_{1}\in X-\{\ell_0,\ell_2,\ell_{3},\ldots,\ell_k\}$ such that $\cC_{k}^-=\langle\ell_0,\ell_{1},\ell_2,\ell_{3},\ldots,\ell_k\rangle$ is a nearly-tight caterpillar ladder of $\cN$. 
\end{enumerate}
\end{lemma}

\begin{proof}
Assume that (i) does not hold. We will show that (ii) holds. For the purposes of the proof, if $k=2$, we set $p_3=q_2$. By Lemma~\ref{l:tree-path-down}, there is a unique tree path $P$ from $p_{3}$ to $p_2$. Let $q_{1}$ denote the child of $p_{3}$ that lies on $P$. Since $P$ does not contain a reticulation, $q_{1}$ is a tree vertex and $q_{1}$ has a child $v_{1}$ that does not lie on $P$. Furthermore, applying  Lemma~\ref{l:tree-path-down} again, it follows that every tree path that starts at $p_3$ traverses $p_2$. Hence $v_{1}$ is a reticulation. Let $\ell_{1}$ be a leaf at the end of a tree path that starts at $v_{1}$, let $f_{1}=(q_{1},v_{1})$, and let $e_{1}=(p_{1},v_{1})$ be the other reticulation arc that is directed into $v_{1}$.

Now, let $\cE'$ be an embedding of a phylogenetic $X$-tree $\cT$ in $\cN$ that uses $\{f_{1},f_2,\ldots,f_k\}$. Then
\begin{equation}\label{eq:tree-one}
\cT|\{\ell_0,\ell_{1},\ell_2,\ldots,\ell_k\}=(\ell_0,\ell_{1},\ell_2,\ldots,\ell_k).
\end{equation}
As $f_k$ is non-essential, there exists an embedding $\cE$ of $\cT$ in $\cN$ that uses $e_k$. Furthermore, by Lemma~\ref{l:embedding-down} and Corollary~\ref{2down}, $\cE$ also uses $\{e_2,e_{3},\ldots,e_{k-1}\}$. Assume that $\cE$ uses $f_{1}$. Then $\cT|\{\ell_0,\ell_{1},\ell_2,\ldots,\ell_k\}$ has $\{\ell_0,\ell_2\}$ as a cherry, a contradiction. We may therefore assume that $\cE$ uses $e_{1}$.

\begin{sublemma}\label{s:shortcut-down}
The arc $f_{1}$ is a shortcut.  
\end{sublemma}

\begin{proof}
Let $x$ be a leaf at the end of a tree path that starts at $p_{1}$. First assume that $x\notin \{\ell_0,\ell_2,\ell_{3},\ldots,\ell_k\}$. By Lemma~\ref{l:tree-path-down}, every tree path that starts at $p_3$ traverses $p_2$. Considering $\cE'$, it now follows that $\cT|\{x,\ell_0,\ell_{1},\ell_2,\ldots,\ell_k\}$ is either $(\ell_0,\ell_{1},\ell_2,\ldots,\ell_k,x)$ or it has a cherry in $$\{\{x,\ell_0\},\{x,\ell_2\},\{x,\ell_{3}\},\ldots,\{x,\ell_k\}\}.$$ On the other hand, by considering $\cE$, it follows that $\cT|\{x,\ell_0,\ell_{1},\ell_2,\ldots,\ell_k\}$ has $\{x,\ell_{1}\}$ as a cherry, a contradiction. Hence, $x\in \{\ell_0,\ell_2,\ell_{3},\ldots,\ell_k\}$. If $x\in \{\ell_2, \ell_3, \ldots, \ell_k\}$, then, by considering $\cE$, we deduce that $\cT|\{\ell_0, \ell_1, \ell_2, \ldots, \ell_k\}$ has a cherry in $\{\{\ell_1, \ell_2\}, \{\ell_1, \ell_3\}, \ldots, \{\ell_1, \ell_k\}\}$, contradicting the restriction of $\cT$ in (\ref{eq:tree-one}). Thus $x=\ell_0$, and so one of $e_{1}$ and $f_{1}$ is a shortcut. Assume that $e_{1}$ is a shortcut. Then there is a directed path from $p_{1}$ to $q_{1}$. In particular, $p_{1}$ is an ancestor of $q_k$. Now, by considering $\cE'$, we again have the restriction of $\cT$ in (\ref{eq:tree-one}), whereas, by considering $\cE$, we have $$\cT|\{\ell_0,\ell_{1},\ell_2,\ldots,\ell_k\}=(\ell_0,\ell_2,\ell_3,\ldots,\ell_k,\ell_{1}),$$ a contradiction. Hence $f_{1}$ is a shortcut.\qed
\end{proof}

We next show that $p_1$ lies on the tree path from $p_2$ to $\ell_0$. Let $x$ be a leaf at the end of a tree path starting at $p_{1}$. Clearly, $x\neq \ell_1$. Furthermore, since there is a directed path from $q_{1}$ to $p_{1}$, it follows that $x\not\in\{\ell_{3},\ell_{4},\ldots,\ell_k\}$. Assume that $x\notin\{\ell_0,\ell_2\}$. Then, by considering $\cE'$ and recalling that $f_1$ is a shortcut, $\cT|\{x,\ell_0,\ell_{1},\ell_2\}$ has a cherry in $\{\{x,\ell_0\},\{x,\ell_2\}\}$, whereas, by considering $\cE$, we have that $\cT|\{x,\ell_0,\ell_{1},\ell_2\}$ has $\{x,\ell_{1}\}$ as a cherry, a contradiction. Hence $x\in\{\ell_0,\ell_2\}$. Continue to consider $\cE$. If $x=\ell_2$, then $p_{1}$ lies on the tree path from $v_2$ to $\ell_2$ and $\cT|\{\ell_{0},\ell_1,\ell_2\}=(\ell_{1},\ell_2,\ell_0)$, contradicting the restriction of $\cT$ in (\ref{eq:tree-one}) when considering $\cE'$. Thus $x=\ell_0$. If $p_{1}$ lies on the subpath of $P$ from $q_{1}$ to $p_2$, then, by considering $\cE$, $\cT|\{\ell_0,\ell_{1},\ell_2\}=(\ell_0,\ell_2,\ell_{1})$, again contradicting the restriction of $\cT$ in (\ref{eq:tree-one}) when considering $\cE'$. It now follows that $p_{1}$ lies on the tree path from $p_2$ to $\ell_0$. In particular, we have that

\begin{sublemma}
$\langle \ell_0, \ell_1, \ell_2, \ldots, \ell_k\rangle$ is a loose caterpillar ladder in $\cN$.
\label{loosesub}
\end{sublemma}

\begin{sublemma}
Let $u$ be the child of $q_{1}$ on $P$. Then $u=p_2$.
\end{sublemma}

\begin{proof}
Towards a contradiction, assume that $u\ne p_2$. Then, as $P$ is a tree path, $u$ has a child $w$ that does not lie on $P$. Furthermore, $w$ is a reticulation because, otherwise, there is a tree path in $\cN$ that starts at $p_3$ and does not traverse $p_2$, contradicting Lemma~\ref{l:tree-path-down}. Let $x$ be a leaf at the end of a tree path starting at $w$. Since $\cN$ is tree-child, $x\notin\{\ell_0,\ell_{1},\ell_2,\ldots,\ell_k\}$. Let $\cE'$ be an embedding of a phylogenetic $X$-tree $\cT$ that uses $\{(u,w),e_{1},f_2,f_{3},\ldots,f_k\}$. Then 
\begin{equation}\label{eq:tree-two-and-a-half}
\cT|\{\ell_0,\ell_{1},\ell_2\}=(\ell_0,\ell_{1},\ell_2) \textnormal{ and}
\end{equation}
\begin{equation}\label{eq:tree-two}
\cT|\{x,\ell_0,\ell_{1},\ell_2\}=(\ell_0,\ell_{1},x,\ell_2).
\end{equation}
Since $f_k$ is non-essential, there exists an embedding $\cE$ of $\cT$ in $\cN$ that uses $e_k$. Moreover, by Lemma~\ref{l:embedding-down} and Corollary~\ref{2down}, $\cE$ also uses $\{e_2,e_{3},\ldots,e_{k-1}\}$. 

We continue to consider $\cE$. First, if $\cE$ uses $f_{1}$, then
$$\cT|\{\ell_0,\ell_{1},\ell_2\}=(\ell_0,\ell_2,\ell_{1}),$$
a contradiction to the restriction of $\cT$ in (\ref{eq:tree-two-and-a-half}). Hence $\cE$ uses $e_{1}$. Second, if $\cE$ uses $(u,w)$, then
$$\cT|\{x,\ell_0,\ell_{1},\ell_2\}=(\ell_0,\ell_{1},\ell_2,x),$$
which gives a contradiction to the restriction of $\cT$ in (\ref{eq:tree-two}) when considering $\cE'$. Thus, if $u'$ denotes the parent of $w$ that is not $u$, then $\cE$ uses $(u',w)$. We next consider two cases relative to $\cE$. First, if $(u', w)$ is a shortcut, or if neither $(u, w)$ nor $(u', w)$ is a shortcut, then $\cT|\{x,\ell_0,\ell_{1},\ell_2\}=(\ell_0,\ell_{1},\ell_2,x)$. Second, if $(u, w)$ is a shortcut and $u'$ does not lie on the  tree path from $p_2$ to $p_{1}$, then $\cT|\{x,\ell_0,\ell_{1},\ell_2\}$ is either $(\ell_0,\ell_{1},\ell_2,x)$ or it has a cherry in $\{\{x,\ell_0\},\{x,\ell_{1}\},\{x,\ell_2\}\}$. Both cases contradict the restriction of $\cT$ in (\ref{eq:tree-two}) when considering $\cE'$; thereby implying that $(u, w)$ is a shortcut and $u'$ lies on the tree path from $p_2$ to $p_{1}$.

Now, let $\cE'_1$ be an embedding of a phylogenetic $X$-tree $\cT'$ that uses $\{f_{1},f_2,\ldots,f_k\}$. Then 
\begin{equation}\label{eq:missing-case}
\cT'|\{x,\ell_0,\ell_{1},\ell_2\}=(\ell_0, x, \ell_1, \ell_2)
\end{equation}
regardless of whether $\cE'_1$ uses $(u,w)$ or $(u',w)$. Since $f_k$ is non-essential, there exists an embedding $\cE_1$ of $\cT'$ in $\cN$ that uses $e_k$. Moreover, by (\ref{loosesub}), Lemma~\ref{l:embedding-down}, and Corollary~\ref{2down}, $\cE_1$ also uses $\{e_{1},e_2,\ldots,e_{k-1}\}$. If $\cE_1$ uses $(u,w)$, then $\cT'|\{x,\ell_0,\ell_{1},\ell_2\}=(\ell_0,\ell_{1},\ell_2,x)$ and, if $\cE_1$ uses $(u',w)$, then $\cT'|\{x,\ell_0,\ell_{1},\ell_2\}=(\ell_0,\ell_{1},x,\ell_2)$. Both cases give a final contradiction to the restriction of $\cT'$ in (\ref{eq:missing-case}). Thus $u=p_2$.\qed
\end{proof}

It now follows that $\cC_{k}^-=\langle\ell_0,\ell_{1},\ell_2,\ldots,\ell_k\rangle$ is a nearly-tight caterpillar ladder of $\cN$. This completes the proof of the lemma. \qed
 \end{proof}
 
At last, we establish Theorem~\ref{t:main}.\\

\noindent {\it Proof of Theorem~\ref{t:main}.}
Suppose that there exists a tight caterpillar ladder $\cC=\langle\ell_0,\ell_1,\ell_2,\ldots,\ell_k\rangle$ of $\cN$. It follows from Lemma~\ref{l:one-direction} that each of $e_1$ and $f_k$, the first and last rungs of $\cC$, is non-essential in $\cN$. For the converse, suppose that $e$ is non-essential. Let $v$ be the reticulation of $\cN$ that $e$ is directed into, and let $u$ and $u'$ be the two parents of $v$. By the contrapositive of Lemma~\ref{l:no-shortcut}, we may assume without loss of generality that $(u',v)$ is a shortcut. Since $\cN$ is tree-child, there exist two distinct elements $\ell_0$ and $\ell_1$ in $X$ such that there is a tree path from $u$ to $\ell_0$ and a tree path from $v$ to $\ell_1$. First, assume that $e=(u,v)$.   Then $\cC_1^+=\langle\ell_0,\ell_1\rangle$ is a nearly-tight caterpillar ladder of $\cN$ with $e_1=e$. Moreover, by Lemma~\ref{l:bottom-up}, $\cC_1^+$ is either a tight caterpillar ladder of $\cN$ or there exists an element $\ell_2$ in $X-\{\ell_0,\ell_1\}$ such that $\cC_2^+=\langle\ell_0,\ell_1,\ell_2\rangle$ is a nearly-tight caterpillar ladder of $\cN$. As $X$ is finite, it now follows that, for some $k\in\{1,2,\ldots,|X|-1\}$, after $k-1$ applications of Lemma~\ref{l:bottom-up}, $\cC_{k}^+=\langle\ell_0,\ell_1,\ell_2,\ldots,\ell_k\rangle$ is a  tight caterpillar ladder of $\cN$ in which $e$ is the first rung. Second, assume that $e=(u',v)$. Then, up to relabeling the leaves, using Lemma~\ref{l:top-down} instead of Lemma~\ref{l:bottom-up} and applying an argument that is  analogous to that applied in the previous case gives the desired result. This completes the proof of the theorem.
\qed

\section{Obtaining Tree-Child Networks with No Non-Essential Arcs}
\label{sec:algorithm}

In this section, we show that it takes polynomial time to compute a tree-child network whose reticulation arcs are all essential and that displays the same set of phylogenetic trees as a given tree-child network. We start with three lemmas  that establish properties of collections of tight caterpillar ladders. The first lemma establishes that if two tight caterpillar ladders have a reticulation or spine vertex in common, then they are essentially the same tight caterpillar ladder. Furthermore, the latter two lemmas show how the two sets of caterpillar ladders of two tree-child networks $\cN$ and $\cN'$ differ if $\cN'=\cN\setminus \{e\}$, where $e$ is the first or last rung of a caterpillar ladder of $\cN$.

Let $\cN$ be a tree-child network on $X$, and let $\cC=\langle\ell_0, \ell_1, \ell_2, \ldots, \ell_k\rangle$ and $\cC'=\langle\ell'_0, \ell'_1, \ell'_2, \ldots, \ell'_{k'}\rangle$ be tight caterpillar ladders of $\cN$. Let $\phi$ and $\phi'$ be maps from the vertex sets of $\cC$ and $\cC'$, respectively, to $\cN$ that satisfy (P1)--(P3). Similar to the convention used throughout the proof of Theorem~\ref{t:main}, we will use the following notation throughout this subsection, where $i\in \{1, 2, \ldots, k\}$ and $j\in \{1, 2, \ldots, k'\}$:
\begin{enumerate}[1.]
\item For each reticulation $v_i$ in $\cC$ (resp.\ $v'_j$ in $\cC'$), the vertex $\phi(v_i)$ (resp.\ $\phi'(v'_j)$) in $\cN$ is denoted by $v_i$ (resp.\ $v'_j$).

\item For each vertex $v_i$ (resp.\ $v'_j$) in $\cN$, the parents of $v_i$ (resp.\ $v'_j$) are denoted by $p_i$ and $q_i$ (resp.\ $p'_j$ and $q'_j$), where there is a directed path from $q_i$ to $p_i$ (resp.\ $q'_i$ to $p'_i$) in $\cN$.
\end{enumerate}

\begin{lemma}
Let $\cN$ be a tree-child network on $X$, and let $\cC$ and $\cC'$ be tight caterpillar ladders of $\cN$. Let $u$ be a vertex of $\cN$. If $u$ is either a reticulation or a spine vertex of $\cC$ and $\cC'$, then all reticulation and spine vertices of $\cC$ and $\cC'$ coincide.
\label{coincide}
\end{lemma}

\begin{proof}
Throughout this proof, let $\cC=\langle\ell_0, \ell_1, \ell_2, \ldots, \ell_k\rangle$, and let $\cC'=\langle\ell'_0, \ell'_1, \ell'_2, \ldots, \ell'_{k'}\rangle$. Note that the spine vertices of a caterpillar ladder are all tree vertices. If $k=1$ and $u$ is a reticulation or a spine vertex of $\cC$ and $\cC'$ in $\cN$, then it is clear that the lemma holds. So we may assume that $k\ge 2$ and, similarly, $k'\ge 2$. Suppose that $u$ is a reticulation or spine vertex of $\cC$ and $\cC'$ in $\cN$. Then, for some $i$ and $j$,  one of $v_i=v'_j$, $p_i=p'_j$, and $q_i=q'_j$ holds. But, regardless of which of these holds, it follows that the other two equalities also hold. In turn, as $\cC$ and $\cC'$ are both tight caterpillar ladders of $\cN$, this implies that, unless $i\in \{1, k\}$ or $j\in \{1, k'\}$, we have $v_{i-1}=v'_{j-1}$, $p_{i-1}=p'_{j-1}$, $q_{i-1}=q'_{j-1}$, and $v_{i+1}=v'_{j+1}$, $p_{i+1}=p'_{j+1}$, $q_{i+1}=q'_{j+1}$. Furthermore, if $i=1$, then, as $\cC$ is a tight caterpillar ladder of $\cN$ and, therefore, satisfies (P3), it is easily seen that $j=1$ and so $v_2=v'_2$, $p_2=p'_2$, and $q_2=q'_2$. Similarly, if $i=k$, then $v_k=v'_k$, $p_k=p'_k$, and $q_k=q'_k$. By repeating this process, for $v_{i-1}$ and then $v_{i+1}$, we eventually deduce that $k=k'$ and, for all $i\in \{1, 2, \ldots, k\}$, we have $v_i=v'_i$, $p_i=p'_1$, and $q_i=q'_i$. This completes the proof of the lemma. \qed
\end{proof}

Following Lemma~\ref{coincide}, two tight caterpillar ladders of a tree-child network are called {\em distinct} if they have no reticulation or spine vertex in common. Furthermore, two tight caterpillars are said to be {\em equivalent} if they have the same set of reticulation and spine vertices.

\begin{lemma}
Let $\cN$ be a tree-child network on $X$, and let $\cC$ and $\cC'$ be two distinct tight caterpillar ladders of $\cN$. If $e$ is the first or last rung of $\cC$, then $\cC'$ is a tight caterpillar ladder of $\cN\setminus\{e\}$.
\label{tightafter}
\end{lemma}

\begin{proof}
Let $e=(u, v)$ be the first or last rung of $\cC$. Since $\cC$ and $\cC'$ are distinct, $e$ is not a reticulation arc of $\cC'$. Now let $P$ be a tree path of $\cN$. Then $P$ is a tree path of $\cN\setminus\{e\}$ unless $P$ traverses $u$ or starts at $v$. If $P$ traverses $u$, then deleting $u$ from $P$ gives a tree path of $\cN\setminus\{e\}$ while, if $P$ starts at $v$, then deleting $v$ from $P$ gives a tree path of $\cN\setminus\{e\}$ starting at the child of $v$ in $\cN$. Since $\cC'$ is a tight caterpillar ladder of $\cN$ whose spine and reticulation vertices in $\cN$ are distinct from those of $\cC$, it now follows that $\cC'$ is a tight caterpillar ladder of $\cN\setminus\{e\}$.\qed
\end{proof}

\begin{lemma}
Let $\cN$ be a tree-child network on $X$, and let $\cC$ be a tight caterpillar ladder of $\cN$. Let $e$ be the first or last rung of $\cC$. If $\cC'$ is a tight caterpillar ladder of $\cN\setminus\{e\}$, then either $\cC'$ or a tight caterpillar ladder equivalent to $\cC'$ is a tight caterpillar ladder of $\cN$.
\label{tightbefore}
\end{lemma}

\begin{proof}
Let $\cC=\langle\ell_0, \ell_1, \ell_2, \ldots, \ell_k\rangle$. Depending on whether $e$ is the first or last rung of $\cC$, after deleting $e$ to obtain $\cN\setminus\{e\}$ either (i) $p_1$ and $v_1$ are suppressed if $e=(p_1, v_1)$ is the first rung of $\cC$ or (ii) $q_k$ and $v_k$ are suppressed if $e=(q_k, v_k)$ is the last rung of $\cC$. Note that, if (ii) applies and $q_k$ is the root of $\cN$, then after deleting $e$, the vertex $q_k$ and its incident arc are also deleted, and $v_k$ is suppressed. We prove the lemma for when (i) holds. The proof of the lemma for when (ii) holds is similar and omitted.

Assume that (i) holds. Let $R$ and $S$ denote the reticulation and spine vertices of $\cC$ in $\cN$. Since $\cN$ is tree-child, the child of $v_1$ in $\cN$ is either a tree vertex or a leaf, and so $q_1$ is not a spine vertex of a tight caterpillar ladder of $\cN\setminus\{e\}$. It now follows by Lemma~\ref{coincide} that no vertex in $(R-\{v_1\})\cup (S-\{p_1\})$ is a reticulation or spine vertex of a tight caterpillar ladder of $\cN\setminus\{e\}$. Thus $\cC'$ is a tight caterpillar ladder of $\cN$ unless the tree path either from $p'_1$ to $\ell'_0$, or from $v'_j$ to $\ell'_j$ for some $j\in \{1, 2, \ldots, k'\}$ in $\cN\setminus \{e\}$ is no longer a tree path in $\cN$. Note that the existence of $e=(p_1, v_1)$ in $\cN$ effects at most one of these tree paths. Let $P'$ denote such a tree path in $\cN\setminus\{e\}$. Then $P'$ traverses $q_1$. Let $u$ denote the child of $p_1$ in $\cN$ that is not $v_1$, and let $w$ denote the the child of $v_1$ in $\cN$. As $\cN$ is tree-child, each of $u$ and $w$ is either a tree vertex or a leaf. Therefore, as $P'$ traverses $q_1$, it also traverses either $u$ or $w$ in $\cN\setminus\{e\}$. If $P'$ traverses $u$, then modify $P'$ by replacing the subpath $q_1,\, p_2,\, u$ with $q_1,\, p_2,\, p_1\, u$ if $k\ge 2$ and replacing the subpath $q_1\, u$ with $q_1,\, p_1,\, u$ if $k=1$. The resulting path is a tree path in $\cN$, and so $\cC'$ is a tight caterpillar ladder of $\cN$. On the other hand, if $P'$ traverses $w$, then replace the subpath of $P'$ from $q_1$ to a leaf $\ell'_{j'}$, where $j'\in \{0, 1, 2, \ldots, k'\}$, with the (unique) tree path from $q_1$ to $\ell_0$. The resulting path is a tree path in $\cN$, and so $\cC''=\langle\ell'_0, \ell'_1, \ldots, \ell'_{j'-1}, \ell_0, \ell_{j'+1}, \ldots, \ell_{k'}\rangle$ is a tight caterpillar ladder of $\cN$ and $\cC''$ is equivalent to $\cC'$. This completes the proof of the lemma.\qed
\end{proof}

\begin{theorem}\label{t:no-name}
Let $\cN$ be a tree-child network on $X$, and let
$$\cN=\cN_0, \cN_1, \cN_2, \ldots, \cN_m=\cN'$$
be a sequence of tree-child networks such that $\cN_{i+1}$ is obtained from $\cN_i$ by deleting a non-essential reticulation arc for all $i\in \{0, 1, \ldots, m-1\}$. If $\cN'$ has no non-essential reticulation arcs, then $\cN$ has exactly $m$ tight caterpillar ladders. Moreover, $\cN'$ can be obtained from $\cN$ by deleting either the first or last rung of each tight caterpillar ladder of $\cN$.
\end{theorem}

\begin{proof}
By Lemmas~\ref{tightafter} and~\ref{tightbefore}, $\cN_{i+1}$ has exactly one less tight caterpillar ladder that $\cN_i$. The theorem now follows from Theorem~\ref{t:main}.\qed
\end{proof}

We next describe an algorithm that finds a tight caterpillar ladder of a tree-child network $\cN$ relative to a given reticulation of $\cN$ if it exists.

\noindent{\sc Find Tight Caterpillar Ladder}\\
\noindent{\bf Input.} A tree-child network $\cN$ and a reticulation $w$ of $\cN$.\\
\noindent{\bf Output.} A tight caterpillar ladder $\langle\ell_0,\ell_1,\ell_2,\ldots,\ell_k\rangle$  of $\cN$ with $v_1=w$ if it exists and, otherwise, a statement saying that no such ladder exists.

\begin{enumerate}[{\bf 1.}]
\item Set $i=1$ and $v_1=w$.
\item \label{leaf1step} Let $\ell_1$ be a leaf at the end of a tree path starting at $v_1$.
\item Let $u$ and $u'$ be the parents of $v_1$.
\item \label{parents1} If there is a directed path from $u$ to $u'$ in $\cN$, then set $p_1=u'$ and $q_1=u$.
\item \label{parents2} If there is a directed path from $u'$ to $u$ in $\cN$, then set $p_1=u$ and $q_1=u'$.
\item Else, go to Step~\ref{step_no}.
\item \label{leaf0step} Let $\ell_0$ be a leaf at the end of a tree path starting at $p_1$.
\item If $i=1$, then do the following:
\begin{enumerate}
\item If  $q_1,p_1,v_1$ is a tree path in $\cN$, then return $\langle\ell_0,\ell_1\rangle$. 
\item If there is a tree path  $q_1,p_{2},p_1$, and there is a reticulation $v_2$ such that $(p_{2},v_{2})$ is a reticulation arc of $\cN$, then do the following:
\begin{enumerate}
\item Increment $i$ by one.
\item \label{leaf2step} Let $\ell_2$ be a leaf at the end of a tree path starting at $v_2$.
\item Let $q_2$ be the parent of $v_2$ that is not $p_2$.
\end{enumerate}
\item Else, go to Step~\ref{step_no}.
\end{enumerate}
\item \label{step} If $i>1$, then do the following:
\begin{enumerate}
\item If  $(q_i,q_{i-1})$ is an arc of $\cN$, then stop and return $\langle\ell_0,\ell_1,\ell_2,\ldots,\ell_i\rangle$.
\item If there is a tree path $q_i,p_{i+1},q_{i-1}$, and there is a reticulation $v_{i+1}$ such that $(p_{i+1},v_{i+1})$ is a reticulation arc of $\cN$, then do the following:
\begin{enumerate}
\item Increment $i$ by one.
\item \label{leafistep} Let $\ell_i$ be a leaf at the end of a tree path starting at $v_i$.
\item Let $q_i$ be the parent of $v_i$ that is not $p_i$.
\item Go to Step~\ref{step}.
\end{enumerate}
\item Else, go to Step~\ref{step_no}.
\end{enumerate}
\item \label{step_no} Return ``There is no tight caterpillar ladder in $\cN$ with $v_1=w.$''
\end{enumerate}

We finish this section by establishing that it takes polynomial time to obtain a tree-child network whose reticulation arcs are all essential and that displays the same set of phylogenetic trees as a given tree-child network.

\begin{theorem}\label{t:recognizing}
Let $\cN$ be a tree-child network on $X$. It takes time $O(|X|^3)$  to obtain a tree-child network on $X$ from $\cN$ that displays the same set of phylogenetic $X$-trees as $\cN$ and has no non-essential arcs.
\end{theorem}

\begin{proof}
To establish the theorem, we first show that {\sc Find Tight Caterpillar Ladder} works correctly and has running time $O(|X|^2)$. We start by bounding the size of the vertex set and arc set of $\cN$. Let $r$ be the number of reticulations of $\cN$. As $\cN$ is tree-child, it follows from Lemma 2.2 in~\cite{mcdiarmid15} that $r\leq |X|-1$. Hence, by  Corollary 2.3 and Lemma 2.1 of the same paper, $\cN$ has at most $4|X|$ vertices and at most $$3r+2|X|-2\leq 3(|X|-1)+2|X|-2=5|X|-5$$ arcs, respectively. 

Let $w$ be a reticulation of $\cN$. By Lemmas~\ref{l:bottom-up} and~\ref{l:top-down}, {\sc Find Tight Caterpillar Ladder} correctly decides whether or not there exists a tight caterpillar ladder of $\cN$ with $v_1=w$. Next, we turn to the running time analysis of {\sc Find Tight Caterpillar Ladder}. We start by noting that applying a breadth-first search to $\cN$ takes time $O(|X|)$ because the running time of such a search is linear in the sum over the size of the vertex set and the size of the arc set of a given input graph~\cite{cormen01}. Consequently, each of Steps \ref{parents1} and \ref{parents2} takes time $O(|X|)$ since the existence of a directed path from a vertex $s$ to another vertex $s'$ in $\cN$ can be decided by applying a  breadth-first search to $\cN$ that starts at $s$. Now let $t$ be a vertex of $\cN$. As the size of the arc set of $\cN$ is $O(|X|)$, it follows that finding a leaf at the end of a tree path that starts at $t$ takes time $O(|X|)$. Hence,  each of Steps~\ref{leaf1step}, \ref{leaf0step}, \ref{leaf2step}, and~\ref{leafistep} takes time $O(|X|)$. Moreover, each remaining step of the algorithm takes constant time. Since {\sc Find Tight Caterpillar Ladder} iterates at most $r$ times through Step~\ref{step}, where $r\leq |X|-1$, it  follows that {\sc Find Tight Caterpillar Ladder} takes time $O(|X|^2)$.

Next, we apply a breadth-first search to $\cN$ and, at each reticulation, call {\sc Find Tight Caterpillar Ladder}. In this process, we delete either the first or last rung, but not both, of each tight-caterpillar ladder that is returned by a call to  {\sc Find Tight Caterpillar Ladder} and suppress the resulting two degree-two vertices. Let $\cN'$ be the resulting tree-child network on $X$. Since the number of vertices in $\cN$ is $O(|X|)$ and each call of {\sc Find Tight Caterpillar Ladder} takes $O(|X|^2)$ time, the aforementioned process takes $O(|X|^3)$ time in total. Furthermore, it follows from Theorem~\ref{t:no-name} that $\cN'$ displays the same set of phylogenetic $X$-trees as $\cN$ and has no non-essential arcs. This completes the proof of the theorem.\qed
\end{proof}

\section{Non-Essential Arcs in Arbitrary Phylogenetic Networks} \label{sec:hardness}

In this section, we show that, in general, deciding if a reticulation arc of a phylogenetic network is non-essential is computationally hard. More specifically, we show that the following problem is $\Pi_2^P$-complete, that is complete for the second level of the polynomial-time hierarchy~\cite{stockmeyer76}.\\

\noindent {\sc Non-Essential-Arc}\\
\noindent{\bf Input.} A phylogenetic network $\cN$ on $X$ and a reticulation arc $e$ of $\cN$.\\
\noindent{\bf Question.} Is $e$ non-essential?\\

\noindent To establish that  {\sc Non-Essential-Arc} is $\Pi_2^P$-complete, we use a reduction from the next problem, which was recently shown to be $\Pi_2^P$-complete~\cite[Theorem 4.3]{doecker19}. \\

\noindent {\sc Display-Set-Containment}\\
\noindent{\bf Input.} Two phylogenetic networks $\cN_1$ and $\cN_2$ on $X$.\\
\noindent{\bf Question.} Is each phylogenetic $X$-tree that is displayed by $\cN_1$ also displayed by $\cN_2$?\\

Before we are in a position to prove   $\Pi_2^P$-completeness of {\sc Non-Essential-Arc}, we need some additional terminology to generalize the notion of $\cN\setminus\{e\}$ to an arbitrary phylogenetic network. Let $G$ be a directed graph. Up to isomorphism, the {\it full simplification} of $G$ is the directed graph obtained from $G$ by applying the following two operations repeatedly until neither operation is applicable:

\begin{enumerate}
\item Suppress each vertex with in-degree one and out-degree one. 
\item For each pair $u$ and $v$ of vertices, delete all but one arc that joins $u$ and $v$.
\end{enumerate}

Now, let $\cN$ be a phylogenetic network on $X$ with root $\rho$, and let $e$ be a reticulation arc of $\cN$.  In what follows, we apply a $3$-step construction to $\cN$ to obtain, up to isomorphism, the phylogenetic network obtained from $\cN$ by deleting $e$. First, obtain $G_1$ from $\cN$ by deleting $e$. Second, obtain the directed acyclic graph $G_2$ from $G_1$ by deleting each vertex and arc of $G_1$ that does not lie on a directed path from $\rho$ to a leaf in $X$. Observe that $G_2$ does not have any unlabeled vertex with out-degree zero and that the set of vertices of $G_2$ with out-degree zero is $X$. Third, obtain $G_3$ by taking the full simplification of $G_2$. By construction, $G_3$ is a phylogenetic network on $X$ unless $\rho$ has out-degree one, in which case, its unique child has out-degree two. If the exception holds, let $G_3$ be the phylogenetic network on $X$ obtained from the full simplification of $G_2$ by deleting $\rho$ and its incident arc. We refer to $G_3$ as $\cN\setminus\{e\}$.

\begin{theorem}\label{t:hard}
{\sc Non-Essential-Arc} is $\Pi_2^P$-complete.
\end{theorem}

\begin{proof}
We first show that {\sc Non-Essential-Arc} is in $\Pi_2^P$. Let $\cN$ be a phylogenetic network on $X$, and let $e$ be a reticulation arc of $\cN$. To decide if $e$ is a non-essential arc in $\cN$, consider an embedding of a phylogenetic $X$-tree $\cT$ in $\cN$. Then use an NP-oracle to decide if $\cT$ is displayed by  $\cN\setminus\{e\}$. Since $\cN$ and $e$ form a no-instance of {\sc Non-Essential-Arc} if $\cT$ is not displayed by $\cN\setminus\{e\}$, it follows that {\sc Non-Essential-Arc} is in $\Pi_2^P$.

To complete the proof, we establish a reduction from the $\Pi_2^P$-complete problem {\sc Display-Set-Containment} to {\sc Non-Essential-Arc}. The reduction has a similar flavor to that used in~\cite[Theorem 4.5]{doecker19}. Throughout the remainder of the proof, we sometimes label internal vertices of a phylogenetic network. The only purpose of these labels is to make references. Indeed, they should not be regarded as genuine labels as those used for the leaves of a phylogenetic network. Now, let $\cN_1$ and $\cN_2$ be two phylogenetic networks on $X$ that form the input to an instance of {\sc Display-Set-Containment} that asks if each phylogenetic $X$-tree that is displayed by $\cN_1$ is also displayed by $\cN_2$. Furthermore, let $n=|X|$. We next construct a phylogenetic network on leaf set $X\cup\{x,y\}$, where $x, y\not\in X$, that, roughly speaking, contains one copy of $\cN_1$ and $n+1$ copies of $\cN_2$. More precisely, let $G$ be the directed acyclic digraph that is obtained from the caterpillar $(w_0,w_1,w_2,\ldots,w_{n+1})$ by applying the following 4-step process:
\begin{enumerate}
\item Replace $w_0$ with $\cN_1$ by identifying $w_0$ with the root of $\cN_1$ and, for each $i\in\{1,2,\ldots,n+1\}$, replace $w_i$ with $\cN_2$ by identifying $w_i$ with the root of $\cN_2$.
\item For each $\ell_i\in X$ with $i\in \{1,2,\ldots,n\}$, identify all leaves labeled $\ell_i$ with a new vertex $u_i$ and add a new arc $(u_i,\ell_i)$.
\item Let $p$ be the parent of $w_0$ and $w_1$. Subdivide the arc $(p,w_0)$ with a new vertex $v_0$, subdivide the arc $(p,w_1)$ with a new vertex $v_1$, add two new vertices $p_x$ and $x$, and add the three new arcs $(v_0,p_x)$, $(v_1,p_x)$, and $(p_x,x)$.
\item Subdivide the arc $(p,v_0)$ with a new vertex $p_y$, add a new vertex $y$, and add the new arc $(p_y,y)$. 
\end{enumerate} 
To complete the construction, let $\cN$ be a phylogenetic network on $X\cup \{x, y\}$ such that $G$ can be obtained from $\cN$ by contracting (reticulation) arcs. That is, $\cN$ is obtained from $G$ by refining each of the vertices $u_1, u_2, \ldots, u_{n}$ so that every non-root and non-leaf vertex has total degree three. In what follows, we use $v_i$ to denote the parent of $w_i$ for each $i\in\{2,3,\ldots,n+1\}$. By construction, note that the parent of $w_0$ is $v_0$, the parent of $w_1$ is $v_1$, and the parent $p_x$ of $x$ is a reticulation. To illustrate, a possibility for $\cN$ is shown in Figure~\ref{fig:construction}. Since the size of $\cN$ is polynomial in the number of arcs of $\cN_1$ and $\cN_2$, it follows that the construction of $\cN$ takes polynomial time.

\begin{figure}[t]
\center
\scalebox{1.3}{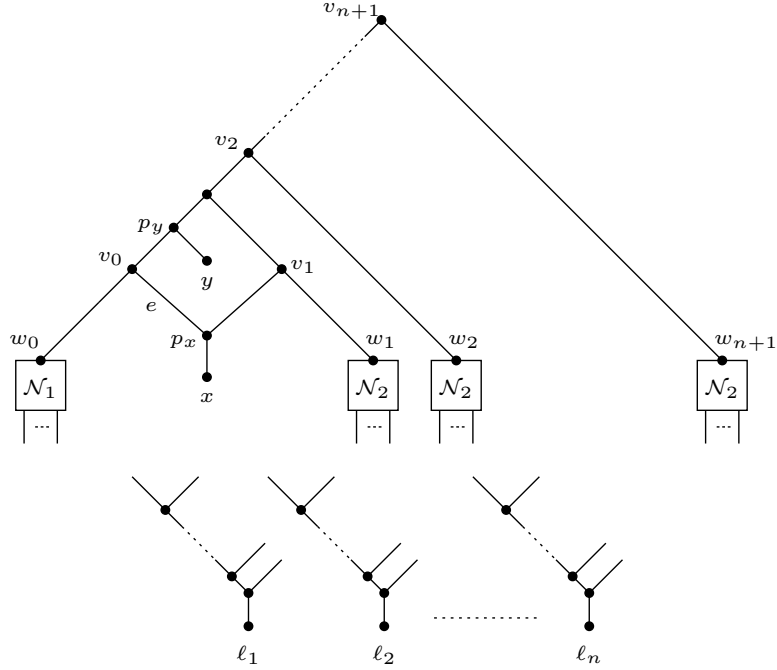}
\caption{The construction of the phylogenetic network $\cN$ on $X\cup \{x,y\}$ as described in the proof of Theorem~\ref{t:hard}, where $X=\{\ell_1,\ell_2,\ldots,\ell_n\}$. The dangling arcs of each square are paired up with the arcs at the bottom part of the figure as described in Step 2 and the subsequent refining procedure of the construction.}
\label{fig:construction}
\end{figure}

\begin{sublemma}
Each  phylogenetic $X$-tree that is displayed by $\cN_1$ is also displayed by $\cN_2$ if and only if $e=(v_0,p_x)$ is non-essential in $\cN$.
\label{sub}
\end{sublemma}

\begin{proof}
First, suppose that not every phylogenetic $X$-tree that is displayed by $\cN_1$ is also displayed by $\cN_2$. Then, there exists a phylogenetic $X$-tree $\cT_1$ that is displayed by $\cN_1$ and not displayed by $\cN_2$. Let $\cT$ be the phylogenetic $(X\cup\{x,y\})$-tree such that $\cT|X=\cT_1$ and, for each pair $\ell,\ell'$ of elements in~$X$, 
\begin{equation}\label{eq:hard-one}
\cT|\{\ell,\ell',x,y\}=(\ell,\ell',x,y). 
\end{equation}
Then there exists an embedding $\cE$ of $\cT$ in $\cN$ that uses $e$ and $(v_0,w_0)$, and no arc in $\{(v_1,w_1),(v_2,w_2),\ldots,(v_{n+1},w_{n+1})\}$. To establish that $e$ is essential, we next show that there exists no embedding of $\cT$ in $\cN\setminus\{e\}$. Note that $\cN\setminus\{e\}$ is the phylogenetic network obtained from $\cN$ by deleting $e$ and suppressing $v_0$ and $p_x$. Assume that there exists an embedding $\cE'$ of $\cT$ in $\cN\setminus\{e\}$. Observe that $\cE'$ uses $(v_1,p_x)$ and $(p_y,y)$. Moreover, as $\cN_2$ does not display $\cT_1$, if follows that $\cE'$ uses at least two distinct arcs $(v_i,w_i)$ and $(v_j,w_j)$ with $i,j\in\{0,1,2,\ldots,n+1\}$. Without loss of generality, we may assume that $i<j$. Let $\ell_i$ (resp. $\ell_j$) be an element of $X$ such that $\ell_i$ (resp. $\ell_j$) is a descendant of $w_i$ (resp. $w_j$) in $\cE'$. By considering $\cE'$, a straightforward check shows that $\cT|\{\ell_i,\ell_j,x,y\}$ is a phylogenetic tree that has a cherry that contains $x$ or $y$, a contradiction to the restriction of $\cT$ in~(\ref{eq:hard-one}). Hence $e$ is essential.

Second, suppose that every phylogenetic $X$-tree that is displayed by $\cN_1$ is also displayed by $\cN_2$. Let $\cT$ be a phylogenetic $(X\cup\{x,y\})$-tree that is displayed by $\cN$ with the property that there exists an embedding $\cE$ of $\cT$ in $\cN$ that uses $e$. To complete the proof, we show that there exists an embedding $\cE'$ of $\cT$ in $\cN\setminus\{e\}$. For each $i\in\{0,1,2,\ldots,n+1\}$, let $X_i$ be the subset of $X$ that contains precisely each element in $X$ that is a descendant of $w_i$ in $\cE$. By the Pigeonhole Principle, there exists at least one $i\in\{1,2,\ldots,n+1\}$ for which $X_i$ is empty. Furthermore $X_0$ may or may not be empty. Let $j$ be the minimum element in $\{1,2,\ldots,n+1\}$ for which $X_j=\emptyset$. For each $i\in\{0,1,2,\ldots,n\}$ with $X_i\ne\emptyset$, we next consider the pendant subtree $\cT|(X_i\cup\{x,y\})$ of $\cT$. Recalling that $\cE$ uses $e$, observe that $\cT|(X_i\cup\{x,y\})$ has cherry $\{x,y\}$ for each $i\in\{1,2,\ldots,n+1\}$. Now, since each phylogenetic $X$-tree that is displayed by $\cN_1$ is also displayed by $\cN_2$, it follows from the construction of $\cN$ that, for each $i\in\{0,1,2,\ldots,n\}$ with $X_i\ne\emptyset$, there is an embedding of the pendant subtree $\cT|(X_i\cup\{x,y\})$  in $\cN$ that uses $(v_{i+1},w_{i+1})$ and $(v_1,p_x)$. Collectively, this implies that, if $X_0\ne\emptyset$, then there exists an embedding of the pendant subtree $$\cT|(X_0\cup X_1\cup X_2\cup\cdots\cup X_{j-1}\cup\{x,y\})$$ in $\cN$ that uses $\{(v_{1},w_{1}),(v_{2},w_{2}),\ldots,(v_{j},w_{j})\}$ and $(v_1,p_x)$. On the other hand, if $X_0=\emptyset$, then there exists an embedding of the pendant subtree $$\cT|(X_1\cup X_2\cup\cdots\cup X_{j-1}\cup\{x,y\})$$ in $\cN$ that uses $\{(v_{2},w_{2}),(v_{3},w_{3}),\ldots,(v_{j},w_{j})\}$ and $(v_1,p_x)$. It is now straightforward to check that there exists an embedding $\cE'$ of $\cT$ in \mbox{$\cN\setminus\{e\}$}.\qed
\end{proof}
\noindent The result (\ref{sub}) completes the proof of Theorem~\ref{t:hard}.\qed
\end{proof}

\bigskip

\noindent{\bf Acknowledgments.} The second author was supported by the New Zealand Marsden Fund.

\end{document}